\documentclass{birkjour}
\usepackage{amscd,amsmath}
\usepackage{amssymb}
\usepackage{amsthm}
\usepackage{color}
\usepackage{hyperref}
\usepackage{enumerate}
\usepackage{geometry}
\theoremstyle{definition}  
\newtheorem{theorem}{Theorem}[section]
\newtheorem{definition}[theorem]{Definition}
\newtheorem{proposition}[theorem]{Proposition}

\newtheorem{corollary}[theorem]{Corollary}
\newtheorem{remark}[theorem]{Remark}
\newtheorem{example}[theorem]{Example}

\numberwithin{equation}{section}

\def\({\left (}
\def\){\right )}
\def\<{\left<}
\def\> { \right>}
\def\e{\eqref}
\def\<{\left<}
\def\>{\right>}
\def\trace{{\rm Trace}}
\def\Ric{{\rm Ric}}

\def\Div{{\rm div}}

\def\grad{grad}

\def\Nkt{{\mathcal M}_{\kappa,\tau}}

\begin{document}

\title[Biconservative submanifolds]{Recent development in biconservative \\submanifolds}

\author[B.-Y. Chen]{Bang-Yen CHEN}
\address{Department of Mathematics, Michigan State University, East Lansing,
Michigan 48824-1027, USA}
\email{chenb@msu.edu}

\begin{abstract}  A submanifold $\phi:M\to \mathbb E^{m}$ is called {\it biharmonic} if it satisfies $\Delta^{2}\phi=0$ identically, according to the author. On the other hand, G.-Y. Jiang studied biharmonic maps between Riemannian manifolds as critical points of the bienergy functional, and proved that biharmonic maps $\varphi$ are characterized by vanishing of bitension $\tau_{2}$ of $\varphi$. During last three decades there has been a growing interest in the theory of biharmonic  submanifolds and biharmonic maps. 

The study of $H$-submanifolds of $\mathbb E^{m}$ were derived from  biharmonic submanifolds by only requiring the vanishing of the tangential component of $\Delta^{2}\phi$.
In 2014, R. Caddeo et. al.  named  a submanifold $M$ in any Riemannian manifold ``biconservative'' if the stress-energy tensor $\hat S_{2}$ of bienergy satisfies ${\rm div}\, \hat S_{2}=0$.  Caddeo et. al.  also shown that a Euclidean submanifolds is an $H$-submanifold if and only if the tangential component of $\tau_{2}$ vanishes and hence the notions of $H$-submanifolds and  of biconservative submanifolds coincide for Euclidean submanifolds.
  The first results on biconservative hypersurfaces were proved by T. Hasanis and  T. Vlachos, where they called such hypersurfaces  {\it H-hypersurfaces} in 1995. 
Since then biconservative submanifolds has attracted many researchers and a lot of interesting results were obtained. 
The aim of this article is to provide a comprehensive survey on recent developments on biconservative submanifolds done most during the last decade. 
\end{abstract} 

\maketitle


\section{Introduction}

As a special topic in the theory of finite type submanifolds, the notion of biharmonic submanifold $M$ of a Euclidean $m$-spaces $\mathbb{E}^{m}$ was introduced in the 1980s by the author as an isometric immersion $\phi: M\to \mathbb{E}^{m}$ of a Riemannian $n$-manifold $M$ into $\mathbb{E}^{m}$ with biharmonic position vector field $\phi$, i.e., $\Delta^{2} \phi=0$  (see, e.g., \cite{chen91,D92}), where $\Delta $ denotes the Laplacian of $M$. According to the famous Beltrami's formula $\Delta \phi=-n H$, the condition ``$\Delta^{2}\phi=0$'' is equivalent to ``$\Delta H=0$'', where $H$ denote the mean curvature vector.

 It is well-known that biharmonic surfaces in $\mathbb{E}^{3}$ are minimal, which leads to  the following conjecture of the author  \cite{chen91} posed in 1991. 
 
\vskip.1in
{\bf  Biharmonic Conjecture:} {\it Biharmonic submanifolds of Euclidean spaces are minimal.}

\vskip.1in
 This conjecture was verified for many cases in a lot of papers under additional
geometric properties (see, e.g., \cite{Chen96,book15,di89,D92,FO22,FHZ21, Fu22,HV95,MOR16,OC20}). However, this conjecture remains open til now.
 
Independently, biharmonic maps were studied by G.-Y. Jiang in \cite{Jiang 2009, Jiang}, in much abstract form; in terms of a variational problem for maps between Riemannian manifolds as critical points of the  $L^{2}$-norm of the tension field.  More precisely,  a biharmonic map $\varphi: M\to N$ between two Riemannian manifolds is a critical point of the {\it bienergy functional}:
\begin{equation*} E_{2}(\varphi)=\frac{1}{2}\int_{M}\|\tau (\varphi)\|^{2}dv,
\end{equation*}
where $M$ is a compact manifold and $\tau (\varphi) = {\rm Trace}(\nabla d\varphi)$ is the tension field of the map $\varphi$.  Jiang proved in \cite{Jiang 2009} that biharmonic maps are characterized by the Euler-Lagrange equation or the biharmonic equation:
\begin{equation}\label{1.1}\tau_{2} (\varphi) = -\Delta^{\phi} \tau (\varphi)-{\rm Trace} \,R^{N}(d\varphi(\cdot), \tau (\varphi))d\varphi(\cdot) = 0,\end{equation}
where $\tau_{2} (\varphi)$ is called the {\it bitension} of $\varphi$, $R^{N}$ is the Riemann curvature tensor of $N$,
and $$\Delta^{\phi}Z={\rm Tr}(\nabla^{\phi}\nabla^{\phi}-\nabla^{\phi}_{\nabla^M})Z.$$

In 2002, R. Caddeo, S. Montaldo and C. Oniciuc  \cite{Cad02} observed that the condition ``$\tau_{2} (\varphi)=0$''  is equivalent to the condition  ``$\Delta^{2} \varphi=0$'' for  Euclidean submanifolds.
Consequently, these two definitions of biharmonicity coincide when the ambient space is Euclidean. During last three decades there has been a growing interest in the theory of biharmonic  submanifolds and biharmonic maps (see, e.g., \cite{OC20}).

The notion of biconservative submanifolds is derived from the theory of biharmonic submanifolds by only requiring the vanishing of the tangential component of $\Delta^{2}\phi$ for Euclidean submanifolds or the tangential component of $\tau_{2} (\varphi)$ for arbitrary Riemannian submanifolds.  
The first results on biconservative hypersurfaces were obtained by T. Hasanis and T. Vlachos in \cite{HV95}, where biconservative hypersurfaces are called {\it $H$-hypersurfaces}.

One of advantages of studying biconservative submanifolds is that one will reveal the influence of the tangential component of the biharmonic equation on submanifolds. Hence, it  may give us informations on the biharmonic conjecture. Moreover, it would provide some informations that how much one can rely on just the vanishing of tangential part of the biharmonic equation as well as how many results remain valid under this half condition. 

Since the last decade, the study of biconservative submanifolds attracts many researchers and a lot of interesting results have been obtained. 
The main purpose of this paper is to present a comprehensive survey of recent developments on biconservative submanifolds. It
is author's intention that this paper would provide useful references for graduate
students as well as for researchers working on this interesting and very active research topic in differential geometry.

\section{Stress-Energy tensor and biconservative submanifolds}\label{S2}

According  to D. Hilbert in \cite{Hilbert1914}, the {\it stress-energy tensor} associated with a variational problem is a symmetric $2$-covariant tensor $\hat S$ conservative at critical points, i.e., ${\rm div}\, \hat S=0$. In the context of harmonic maps $\varphi$  between two Riemannian manifolds the stress-energy tensor have been studied in details by P. Baird and J. Eells in \cite{BE81} and A. Sanini in \cite{AS1983}.

The Euler-Lagrange equation associated to the energy is equivalent to the vanishing of the tension field $\tau_1(\varphi)={\rm Trace}\,(\nabla d\varphi)$ (see, e.g., \cite{ES64}), and the tensor
$$\hat S=\frac{1}{2}\vert d\varphi\vert^2 g - \varphi^{\ast}\tilde g,$$
which satisfies ${\rm div}\, \hat S=-\langle\tau_1(\varphi), d\varphi\rangle$. Hence, we have ${\rm div}\, \hat S=0$ when the map $\varphi$ is harmonic.

For isometric immersions, the condition ${\rm div}\, \hat S=0$ is always satisfied, because $\tau_1(\varphi)$ is normal.
The study of the stress-energy tensor for the
bienergy was initiated G.-Y. Jiang \cite{Jiang87} and afterwards  by  E. Loubeau, S. Montaldo, and C. Oniciuc  \cite{LMO08}.
Its expression is given by
\begin{eqnarray*}\hat S_2(X,Y)&=&\frac{1}{2}\vert\tau_1(\varphi)\vert^2\langle X,Y\rangle+\langle d\varphi,\nabla\tau_1(\varphi)\rangle \langle X,Y\rangle \\
\nonumber && -\langle d\varphi(X), \nabla_Y\tau_1(\varphi)\rangle-\langle d\varphi(Y), \nabla_X\tau_1(\varphi)\rangle, \end{eqnarray*}
and it satisfies the condition
\begin{equation}\label{2.1} {\rm div}\,\hat S_2=-\langle\tau_2(\varphi),d\varphi\rangle, \end{equation}
and thus it conform to the principle of a  stress-energy tensor for the bienergy.

For an isometric immersion  $\varphi:(M,g)\to (\tilde M,\tilde g)$,  equation \eqref{2.1} becomes
\begin{equation} \label{2.2} {\rm div}\,\hat S_2=-\tau_2(\varphi)^{T}.\end{equation}
Based on these facts, R. Caddeo, S. Montaldo, C. Oniciuc and P. Piu  \cite{CMOP14} defined the  notion of {\it biconservative submanifolds} as follows.

\begin{definition} A submanifold of a Riemannian manifold is called {\it biconservative} if it satisfies ${\rm div}\,\hat S_2=0$.
\end{definition}

\section{Basics of submanifolds}\label{S3}

For general references on submanifolds we refer to \cite{book11,book15,book17,book19}.

Let $(M,g)$ be a pseudo-Riemannian $n$-manifold. We denote by $\nabla$ the Levi--Civita connection of $(M,g)$. 
 Then the Riemann curvature tensor $R$ of $M$ is the $(1,3)$ tensor field given by
$$R(X,Y)Z=\nabla_{X}\nabla_{Y} Z- \nabla_{Y}\nabla_{X} Z-\nabla_{[X,Y]} Z$$
for $X,Y,Z\in \mathfrak X(TM)$, where $ {\mathfrak X}(TM)$ denotes the space consists of all (smooth) vector fields in the tangent bundle $TM$ and let $[\;\,, \;]$ be the Lie bracket.
At any point $p\in M$,   for a given basis $\{X,Y\}$ of a plane section $\pi \subset T_{p}M$, we put
$$Q(X,Y)=g(X,X)g(Y,Y)-g(X,Y)^{2}.$$
We call the plane section $\pi$  {\it nondegenerate} whenever $Q(X,Y)\ne 0$. 
For a nondegenerate plane section $\pi$, the~sectional curvature $K(\pi)$ is defined by
$$K(\pi)=\frac{g(R(X,Y)Y,X)}{Q(X,Y)}.$$
The {\it Ricci tensor} $Ric_{p}$ at a point $p\in M$ is given by
 $$\Ric_{p}(X,Y)=\sum_{i=1}^n\epsilon_{i }g(R(E_i,X)Y,E_i),$$
  where $\{E_1,\dots,E_n\}$ is an orthonormal basis of $T_{p}M$;
 and the {\it scalar curvature}  is defined by 
 $$\tau=\sum_{i<j}K(E_{i}, E_{j}).$$
 A Riemannian manifold $(M,g)$ with $\dim M\geq 3$ is called an {\it Einstein manifold} if its Ricci tensor satisfies $\Ric=c g$ for some constant $c$. 
 
\begin{remark} In this article, by a {\it real space form} we mean a Riemannian manifold of constant sectional curvature. We  denote an $m$-dimensional real space form of constant curvature $c$ simply by $R^{m}(c)$.
\end{remark}

For an isometric immersion $\phi: M\to \tilde M$ of a  pseudo-Riemannian submanifold $(M,g)$ into  another  pseudo-Riemannian manifold $(\tilde M,\tilde g)$,  we
denote  the Levi--Civita connections of $M$ and $\tilde M$ by $\nabla$ and $\widetilde\nabla$, respectively. Then the~formulas of Gauss and Weingarten of $M$ in $\tilde M$ are given respectively by
\begin{equation} \label{3.1} \tilde \nabla_XY=\nabla_XY+h(X,Y),\end{equation}
\begin{equation} \label{3.2}  \widetilde\nabla_X\xi=-A_\xi X+D_X \xi,\end{equation}
 where  $X,Y\in \mathfrak X(TM)$ and $\xi\in \mathfrak X(T^{\perp}M)$, where $\mathfrak X(T^{\perp}M)$ is the space consists of all (smooth) normal vector fields of $M$. 
 
 The $h,A$ and $D$ given in \e{3.1}--\e{3.2} are called the {\it second fundamental form},  {\it shape operator} and {\it normal connection} of $M$, respectively.
The~shape operator $A_{\xi}$ is a self-adjoint endomorphism of $T_pM$. Furthermore, the shape operator and the second fundamental form are related by $$ \<h(X,Y),\xi\>=\<A_{\xi}X,Y\>,$$ where $\<\,\;,\;\>$ denotes the inner product associated with the metric.
The {\it mean curvature vector} of $M$ is given by
\begin{equation}  H=\frac{1}{n}\,{\rm Trace}\, h =\frac{1}{n}\sum_{i=1}^{n}\epsilon_{j}h(E_{j},E_{j}).\end{equation}

A submanifold $M$ is called {\it CMC} if $M$ has constant mean curvature $\alpha$, where  
$\alpha=\|H\|$, and $M$ is called  {\it PMCV\/} it its mean curvature vector $H$ is parallel in the normal bundle, i.e., $DH=0$ (see, e.g., \cite{Chen72,Chen10}).
Moreover,  $M$ is called {\it minimal} (resp., {\it totally geodesic}) if its mean curvature vector (resp., its second fundamental form) vanishes identically. Furthermore, $M$ is called {\it totally umbilical} if it  satisfies
$$h(X,Y)=g(X,Y)H,\;\; \forall X,Y\in \mathfrak X(TM);$$
and  $M$ is called {\it pseudo-umbilical} if it shape operator satisfies $$A_{H}=\lambda I$$ for some function $\lambda\in C^{\infty}(M)$, where $I$ is the identity map.

\vskip.05in
Let $R$ and $\tilde R$ denote  the Riemann curvature tensors of $M$ and of $\tilde M$, respectively, and $R^D$ the curvature tensor of  the normal bundle $T^\perp M$ defined by
\begin{equation*}  R^D(X,Y)\xi=D_X D_Y \xi-D_Y D_X \xi-D_{[X,Y]}\xi. \end{equation*}
Then the equations of Gauss, Codazzi and Ricci of $M$ are given respectively by
\begin{equation}\begin{aligned}\label{3.4} & R(X,Y;Z,W)= {\tilde R}(X,Y;Z,W)+\tilde g(h(X,W),h(Y,Z))
\\& \hskip1.1in -\tilde g(h(X,Z),h(Y,W)),\end{aligned}\end{equation}
\begin{equation} \label{3.5}({\tilde R}(X,Y)Z)^\perp= (\bar\nabla_X h)(Y,Z)-(\bar\nabla_Y h)(X,Z),\end{equation}
\begin{equation} \label{3.6} R^D(X,Y;\xi,\eta)= {\tilde R}(X,Y;\xi,\eta)+g([A_\xi,A_\eta ](X),Y), \end{equation}
where  $[A_\xi, A_\eta]=A_\xi A_\eta- A_\eta A_\xi$ and $\overline\nabla h$ is defined by
$$(\overline\nabla_X h)(Y,Z) = D_X h(Y,Z) - h(\nabla_X Y,Z) - h(Y,\nabla_X Z).$$
In particular, if the ambient space $\tilde M$ is a real space form $R^{m}(c)$, then the equations of Gauss, Codazzi and Ricci reduce to
\begin{equation}\begin{aligned} \label{3.7} &g(R(X,Y)Z,W) =  g(A_{h(Y,Z)} X,W) -  g(A_{h(X,Z)}Y,W)
\\& \hskip 1.1in +c\,(g(X,W)g(Y,Z)-g(X,Z)g(Y,W)),\end{aligned}\end{equation}
\begin{equation} \label{3.8}(\overline\nabla_X h)(Y,Z) = (\overline\nabla_Y h)(X,Z),\end{equation}
\begin{equation}\label{3.9} \tilde g( R^{\perp}(X,Y)\xi,\eta ) = g([A_{\xi},A_{\eta}]X,Y).
\end{equation}

\vskip.05in
Let $(\mathbb E^m_s,g_{0})$ denote the pseudo-Euclidean $m$-space whose metric tensor is given by
\begin{equation} g_0=-\sum_{i=1}^s dx_i^2+\sum_{j=s+1}^n dx_j^2,\end{equation}
where $\{x_1,\ldots,x_m\}$ is the rectangular coordinate system of
$\mathbb E^m_s$. Then $(\mathbb E^m_s,g_0)$ is called the {\it pseudo-Euclidean $m$-space} with
index $s$.

For $c>0$, if we put
\begin{equation}  S^m_s(c)=\left\{x\in \mathbb E^{m+1}_s: \<x,x\>=\frac{1}{c^{2}}>0\right\},\end{equation}
\begin{equation} H^m_s(-c)=\left\{x\in \mathbb E^{m+1}_{s+1}: \<x,x\>= -\frac{1}{c^{2}}<0\right\},\end{equation}
where $\<\;,\:\>$ is the indefinite inner product on $\mathbb E^m_s$, then $S^m_s(c)$ and $H^m_s(-c)$
are complete pseudo-Riemannian manifolds with index $s$ and of constant
sectional curvature $c>0$ and $-c<0$, respectively. 
In particular the pseudo-Riemannian manifolds $\mathbb E^m_1,S^m_1(c)$ and $H^m_1(-c)$ are called  the {\it Minkowski spacetime},  {\it de Sitter spacetime} and {\it anti-de Sitter
spacetime}, respectively. These three spacetimes are known as the {\it Lorentzian space forms.}

In the following, by a {\it nondegenerate} submanifold of a Lorentzian space form, we mean a submanifold which is spacelike or timelike.

\section{Two basic formulas and biconservative equations}\label{S4}

The following {formula of Beltrami} is basic in submanifold theory. 

\begin{theorem}\label{T:4.1} {\rm  \cite{beltrami}} {\it If $\phi :M \to \mathbb E^m_s$ is an isometric immersion of a pseudo-Riemannian $n$-manifold $M$ into a pseudo-Euclidean $m$-space $\mathbb E^m_s$, then we have}
\begin{equation}\label{4.1} \Delta \phi=-nH.\end{equation} \end{theorem}

An immediate consequence of Proposition \ref{T:4.1} is the following.

\begin{corollary}\label{C:4.2} {\it Every  spacelike minimal submanifold $M$ in a pseudo-Euclidean space $\mathbb E^m_s$ is non-closed.}
\end{corollary} 
Another  application of Theorem  \ref{T:4.1} is the following.

\begin{corollary} \label{C:4.3} {\it Every spacelike minimal submanifold $M$ in a pseudo-hyperbolic  space $H^m_s(-1)$ is non-closed.} \end{corollary} 

The theory of finite type submanifolds was initiated by the author in late 1990s (see, e.g., \cite{chen79,chen83,chen86,CP87,Chen96,book15}. In particular, the following 
result  is basic and very useful in the study of finite type submanifolds.
 
 \begin{theorem}\label{T:4.4} {\rm \cite{chen79,chen86}} {\it If $\phi :M\to \mathbb E^m_s$ is an isometric immersion of a pseudo-Riemannian $n$-manifold $M$ into $\mathbb E^m_s$, then we have 
\begin{equation}\label{4.2} \Delta H=\Delta^{D}H + \sum_{i=1}^{n}\epsilon_i  h(A_{H}e_{i},e_{i}) + (\Delta H)^T ,\end{equation} where 
\begin{equation}\label{4.3}  (\Delta H)^T=\frac{n}{2} \nabla \| H\|^{2} + 2\, {\rm Trace} (A_{DH})\end{equation}
is the tangential component of $\Delta H$, $\Delta^D$ is the Laplacian associated with the normal connection $D$,  $ \| H\|^{2} =\<H,H\>$, and 
$$ {\rm Trace} (A_{DH})= \sum_{i=1}^{n}  \epsilon_i  A_{D_{E_iH}}E_i$$
for an orthonormal frame  $\{E_1,\ldots,E_n\}$  of $TM$.}
\end{theorem}

Theorem \ref{T:4.4} implies immediately the following. 
 
 \begin{corollary}\label{C:4.5} {\it Let $M$ be a pseudo-Riemannian submanifold of $\mathbb E^m_s$. If $M$ is PMCV,  then it is biconservative.}
 \end{corollary}

For submanifolds in an arbitrary Riemannian manifold,  we have  the following result from \cite{Cad02,Ou10}.

\begin{theorem}\label{T:4.6} {\it If $\phi : M\to \tilde M^{m}$ is an isometric immersion of a Riemannian $n$-manifold $(M,g)$ into another Riemannian manifold $(\tilde M^{m},\tilde g)$, then
$\phi$ is biharmonic if and only if the following two equations are satisfied:
\begin{equation}\label{4.4} 2\,{\rm Trace}(\nabla A_{D H})+ \frac{n}{2} \nabla {\| H\|}^2 = -2\sum_{i=1}^n ({\tilde R}(e_i,H)e_i)^{T}, \end{equation}
\begin{equation} \label{4.5}
\Delta^D H-\sum_{i=1}^n h(A_H E_i,E_i)=  \sum_{i=1}^n (  {\tilde R}(E_i, H)E_i)^\perp,
  \end{equation}
where $\{E_1,\ldots,E_n\}$ is an orthonormal local frame of $TM$.}
\end{theorem}

Theorem \ref{T:4.6} implies the following.
 
 \begin{corollary}\label{C:4.7} {\it If $\phi : M\to \tilde M^{m}$ is an isometric immersion of a Riemannian $n$-manifold $(M,g)$ into another Riemannian manifold $(\tilde M^{m},\tilde g)$, then
$\phi$ is biconservative if and only if equation \e{4.4} holds identically.}
\end{corollary} 

Using the Codazzi equation, one easily find the next result.

\begin{proposition}\label{P:4.8} {\it Let $\phi : M^{n}\to \tilde M^{m}$ be a nondegerate submanifold of $\tilde M$. Then we have}
\begin{equation}\notag
\label{subvarietateoarecare1}
\trace (\nabla A_H) = \frac{n}{2}\grad \left( \|H\|^2 \right) + \trace\, A_{\nabla_{\cdot}^\perp H}(\cdot)+ \trace\left(R^N(\,\cdot\,,H)\,\cdot\, \right)^T.
\end{equation}
\end{proposition}

Because every biharmonic submanifold is biconservative, the family of biconservative submanifolds
is much richer than the family of biharmonic submanifolds. 

When the ambient space $\tilde M^{m}$ is a real space form, Corollary \ref{C:4.7} reduces to the following.

 \begin{corollary}\label{C:4.9} {\rm \cite{CMOP14}} {\it If $M$ is an $n$-dimensional submanifold of a real space form of constant curvature $c$, then $M$ is biconservative if and only if we have}
 \begin{equation}\label{4.6} 2\,{\rm Trace}(\nabla A_{D H})= -\frac{n}{2} \nabla{\| H\|}^2\end{equation}
\end{corollary} 

In particular, if $M$ is a hypersurface, Corollary \ref{C:4.9} reduces to

 \begin{corollary}\label{C:4.10} {\rm \cite{CMOP14}} {\it If $M$ is a hypersurface of a real space form $R^{n+1}(c)$ of constant curvature $c$, then $M$ is biconservative if and only if we have}
 \begin{equation}\label{4.7.1} A(\nabla \alpha)= -\frac{n}{2} \alpha \nabla \alpha,\end{equation}
where $\alpha=\| H\|$.
\end{corollary}

Analogous to Corollary \ref{C:4.5},  Corollary \ref{C:4.9} implies the following.

 \begin{corollary}\label{C:4.11} {\it Let $M$ be a submanifold of a real space form. If $M$ has parallel mean curvature vector, then $M$ is biconservative.}
 \end{corollary}

  Clearly, Corollary \ref{C:4.11} implies 
 
 \begin{corollary}\label{C:4.12} {\it CMC hypersurfaces of a real space form are biconservative.}
 \end{corollary}

\section{Geometry of biconservative surfaces in $R^{3}(c)$}\label{S5}

Let $(M^{2},g)$ be a Riemannian 2-manifold. Then we can write the metric $g$ in terms of an isothermal coordinates as $g=e^{2\rho}(dx^{2}+dy^{2})$, where $\rho\in C^{\infty}(M^{2})$. As usual, let us put 
$$\frac{\partial}{\partial {z}}=\frac{1}{2}\(\frac{\partial}{\partial {x}}-i\frac{\partial}{\partial {y}}\), \;\;  \;\frac{\partial}{\partial {\bar z}}=\frac{1}{2}\(\frac{\partial}{\partial {x}}+i \frac{\partial}{\partial {y}}\) .$$
 Then $g( A_{H}(\tfrac {\partial}{\partial {z}}),\tfrac{\partial}{\partial {z}})$ is called {\it holomorphic} if  
 $\tfrac{\partial}{\partial {\bar z}} g\!\( A_{H}\! \(\tfrac{\partial}{\partial {z}}\),\tfrac{\partial}{\partial {z}}\)=0$ holds.

\subsection{Characterizations of biconservative surface}\label{S5.1}

For  surfaces of an arbitrary Riemannian manifold, we have the following characterization result from {\rm \cite{MOR16.2,Ni17}}.

\begin{theorem}\label{T:5.1} {\it  Let $\phi: (M^{2},g)\to (\tilde M^{m},\tilde g)$ be a surface in a Riemannian $m$-manifold $\tilde M^{m}$. Then the following four statements are equivalent:}
\begin{enumerate}
 \item[{\rm (a)}] {\it $M$ is a biconservative surface;}

\item[{\rm (b)}] {\it $g( A_{H}(\tfrac {\partial}{\partial {z}}),\tfrac{\partial}{\partial {z}})$ is holomorphic;}
\item[{\rm (c)}] {\it $A_{H}$ is a Codazzi tensor.}
\item[{\rm (d)}] {\it the mean curvature $\alpha=\|H\|$ is constant.}
\end{enumerate}
\end{theorem}

\begin{remark} In the case of biharmonic surfaces, the fact that $g( A_{H}(\tfrac {\partial}{\partial {z}}),\tfrac{\partial}{\partial {z}})$ 
is holomorphic if and only if $M$ is CMC  was proved in \cite{LO14}.
\end{remark}

From Corollary \ref{C:4.12} we know that CMC surfaces in any real 3-space form are biconservative, thus the study of non-CMC biconservative surfaces is more interesting. 

\begin{definition}
A maximal biconservative surfaces with $\nabla \alpha\neq 0$ everywhere is called a {\it standard biconservative surface}, and the domain of their defining immersion endowed with the induced metric is called an {\it abstract standard biconservative surface}.
\end{definition}

\subsection{Biconservative surfaces in  $\mathbb E^{3}$}\label{S5.2}

In \cite{Ni16}, S. Nistor constructed some complete biconservative surfaces in $\mathbb E^{3}$ as follows.

\begin{theorem} \label{T:5.4} {\it If $(\mathbb{R}^2,g=c (\cosh u)^6(du^2+dv^2))\, (c>0)$ is a surface, then we have:}
\begin{itemize}
\item[(a)] {\it the metric $g$ on $\mathbb{R}^2$ is complete;}
\item[(b)] {\it the Gauss curvature  $G=G(u)$ of the surface satisfies
           $$
            G(u)=-\frac{3}{c\left(\cosh u\right)^8}<0, \;\;  G^\prime(u)=\frac{24 \sinh u}{c\left(\cosh u\right)^9},
           $$
           and therefore $\nabla G\neq 0$ at any point of $\mathbb{R}^2\setminus Ov$;}

\item[(c)] {\it the immersion $\phi:(\mathbb{R}^2,g)\to \mathbb{R}^3$ given by
    $$
    \phi(u,v)=\left(f(u)\cos 3v, f(u)\sin 3v, k(u)\right)
    $$
        is a biconservative surface in $\mathbb{R}^3$, where}
    $$
    f(u)=\frac{\sqrt{c}}{3}\left(\cosh u\right)^3, \;\;
    k(u)=\frac{\sqrt{c}}{2}\left(\frac{1}{2}\sinh 2u+u\right), \;\;  u \in \mathbb{R}.
    $$
\end{itemize}
\end{theorem}

 In \cite{Ni16}, S. Nistor determined (locally) Riemannian surfaces which admit biconservative immersions with $\nabla \alpha\ne 0$ in $\mathbb E^{3}$ ($\alpha=\| H\|$). She also classified complete biconservative surfaces in  $\mathbb E^{3}$ as the next two results.
 
\begin{theorem} \label{T:5.5} {\it Let $({\mathbb R}^{2},g_{c}=c (\cosh u)^{6}(du^{2}+dv^{2}))$ be a Riemannian surface, where $c$ is a positive constant. Then we have}
\begin{enumerate}
 \item[{\rm (a)}] {\it the metric on ${\mathbb R}^{2}$ is complete;}
 \item[{\rm (b)}] {\it  $\nabla \alpha\ne 0$ at any point of ${\mathbb R}^{2}\setminus Ov$, where $Ov$ is the $v$-axis in ${\mathbb R}^{2}$;}
 \item[{\rm (c)}]  {\it  the immersion $\phi_{c}:({\mathbb R}^{2},g_{c})\to \mathbb E^{3}$ given by
 $$\phi_{c}:(u,v)=(f(u)\cos 3v, f(u)\sin 3v, k(u))$$
is biconservative in  $\mathbb E^{3}$, where
 $$f(u)=\frac{c^{1/2}}{3}(\cosh u)^{3},\;\; k(u)=\frac{c^{1/2}}{2}\(\frac{1}{2} \sinh 2u +u\),\;\; u\in \mathbb R.$$
 } \end{enumerate}
\end{theorem}

\begin{theorem} \label{T:5.6}  {\it Any two complete biconservative surfaces in $\mathbb E^{3}$ differ by a homothety of $\mathbb E^{3}$.} \end{theorem}

For  biconservative surfaces of revolution in $\mathbb E^{3}$, R. Caddeo, S. Montaldo, C. Oniciuc and P. Piu proved the following result.

\begin{theorem}\label{T:5.7} {\rm \cite{CMOP14}} {\it
Let $M$ be a biconservative surface of revolution in ${\mathbb E}^3$ with non-constant mean curvature. Then, locally, the surface can be parametrized by
$\phi_c(u,v)=\big(u \cos v, u \sin v, f(u)\big),$
where 
$$f(u)=\frac{3}{2 c}\(u^{1/3} \sqrt{cu^{2/3}-1}+\frac{1}{\sqrt{c}} \ln \left[2(c u^{1/3} + \sqrt{c^2 u ^{2/3}-c})\right]\)\,,$$
with $c$  a positive constant and $u\in(c^{-3/2},\infty)$. The parametrization $\phi_c$ consists of a family of biconservative surfaces of revolution any two of which  are not locally isometric.}
\end{theorem}
 
A {\it regular surface} in $\mathbb E^{3}$ is a  surface which is defined by an embedding. 
 S. Nistor and C. Oniciuc  \cite{NO19} obtained the following.

\begin{theorem} \label{T:5.8} {\it The only complete biconservative regular
surfaces in the Euclidean 3-space $\mathbb E^{3}$ are either CMC surfaces or certain surfaces of revolution. In particular,
any compact biconservative regular surface in $\mathbb E^{3}$ is a round sphere.}
\end{theorem}

\subsection{Biconservative surfaces in  $S^{3}$}\label{S5.3}

For biconservative surfaces in $S^{3}(1)$, we have the next result of S. Nistor.

\begin{theorem}\label{T:5.9} {\rm  \cite{Ni16}}
{\it If $M^2$ is a biconservative surface in ${S}^3(1)\subset \mathbb{R}^4$ with $\alpha>0$ and  $(\nabla \alpha)\neq 0$ everywhere on $M^{2}$ with $\alpha=\| H\|$,  then,  locally, the surface can be parametrized by
\begin{equation*}\label{eq:YC1tilda}
\phi_{c}(u,v)=\gamma(u)+\frac{4\kappa(u)^{-3/4}}{3\sqrt{c}}\left( v_{1}(\cos v -1)+v_{2} \sin v\right),
  \end{equation*}
where $c$ is a real number,  $v_{1}, v_{2}$ are two  orthonormal constant vectors in $\mathbb{R}^4$, $\gamma$ is a curve parametrized by arc-length  satisfying
\begin{equation*}
\label{eq:sigma_prod_scal}
 \langle \gamma(u),v_{1 }\rangle = \frac{4\kappa(u)^{-3/4}}{3\sqrt{c_1}}, \;\; \langle \gamma(u),v_{2}\rangle=0,
\end{equation*}
and, the curvature $\kappa=\kappa(u)$of $\gamma$ is a positive non-constant solution of \begin{equation}\label{k''k}
\kappa^{\prime\prime}\kappa=\frac{7}{4}\left(\kappa^\prime\right)^2+\frac{4}{3}\kappa^2-4\kappa^4.
\end{equation}}
\end{theorem}

By using a methods developed in \cite{Fu15}, and  a similar discussion as  given in \cite{CMOP14}, Y. Fu \cite{Fu15} established explicit classifications of biconservative surfaces in $S^{3}(1)$ as follows.

\begin{theorem}\label{T:5.10} {\it If $\phi : M\to S^{3}(1)\subset \mathbb E^{4}$ is a biconservative surface in $S^{3}(1)$, then $M$ is either a CMC surface or locally a rotational surface defined by
$$\phi(u,v)= \(u \cos v, u \sin v, \sqrt{1-u^{2}}\cos f, \sqrt{1-u^{2}}\sin f\),$$
where $u\in (0,1)$ and}
$$f=\pm \int \frac{3du}{u^{1/3}(1-u^{2})\sqrt{1-9 u^{-2/3}-u^{2}}}.$$
\end{theorem}

\subsection{Biconservative surfaces  in $H^{3}$}\label{S5.4}

 S.  Nistor and C. Oniciuc \cite{NO20} constructed complete, simply connected, non-CMC biconservative surfaces in $H^{3}(-1)$, in an intrinsic and extrinsic way. They constructed
three families of such surfaces. For each surface, the set of points where
the gradient of the mean curvature function does not vanish is dense and has
two connected components. 
In their intrinsic approach, they first constructed a complete simply-connected, abstract surface and proved that it admits a unique biconservative immersion in $H^{3(-1)}$. 

More precisely,  Nistor and Oniciuc \cite{NO20} constructed the family of abstract complete surfaces $(\mathbb R^{2},\tilde g_{-1,c})$, where $c$ is a real constant. In fact, they proved the following (cf. Theorem 4.1 of \cite{NO22}).

\begin{theorem}\label{T:5.11} {\rm  \cite{NO20}}  {\it There exists a unique biconservative immersion}
\begin{equation}\Phi_{}: (\mathbb R^{2},g_{-1,c})\to H^{3}(-1).\end{equation}
\end{theorem}

In  \cite{Fu15}, Fu also established ``explicit'' classifications of biconservative surfaces in $H^{3}(-1)$ as follows.

\begin{theorem}\label{T:5.12} {\it If $\phi : M\to H^{3}(-1)\subset \mathbb E^{4}_{1}$ is a nondegenerate biconservative surface, then $M$ is either CMC or locally given by one of the following surfaces:}
\begin{enumerate}
 \item[{\rm (1}] {\it a rotational surface defined by
 $$ \phi(u,v)= \(\sqrt{1+u^{2}} \cosh f, u \cos v, \sqrt{1+u^{2}} \sinh f,u \sin v, \),$$
 where $u\in (0,\infty)$ and }
 $$f=\pm \int \frac{3du}{u^{1/3}(1+u^{2})\sqrt{1-9 u^{-2/3}+u^{2}}}.$$
 
 \item[{\rm (2)}] {\it a rotational surface defined by
 $$ \phi(u,v)= \(u \cosh v, u\sinh v,\sqrt{u^{2}-1} \cos f, ,\sqrt{u^{2}-1} \sin f\),$$
 where $u\in (1,\infty)$ and }
 $$f=\pm \int \frac{3du}{u^{1/3}(1-u^{2})\sqrt{9 u^{-2/3}+u^{2}-1}}.$$
 
 \item[{\rm (3)}]  {\it a rotational surface defined by
 $$ \phi(u,v)= \(\frac{1}{2}\left\{u(v^{2}+ f^{2})+\frac{1}{u}+u\right\}\! ,\frac{1}{2}\left\{u(v^{2}+ f^{2})+\frac{1}{u}-u\right\}\! ,uv,vf \),$$
 where $u\in (3^{3/4},\infty)$ and }
 $$f= \int \frac{3du}{u^{2}\sqrt{ u^{8/3}-9}}.$$
 \end{enumerate}
\end{theorem}

\subsection{Biconservative surfaces in  $R^{3}(c)$}\label{S5.5}

The next result of Caddeo, Montaldo, Oniciuc and Piu \cite{CMOP14} provided an intrinsic properties of biconservative surfaces satisfying $\nabla\alpha \ne 0$ in $R^{3}(c)$.

\begin{theorem}\label{T:5.13} {\it  If $\phi:M^2\to R^3(c)$ is a biconservative surface satisfying $\nabla \alpha \neq 0$ $(\alpha=\| H\|)$, then the Gauss curvature  $G$ satisfies}
\begin{itemize}
\item [(1)] $G=\det A+c=-3\alpha^2+c,\;\;$;

\item [(2)] {\it $c-G>0$, $\nabla G\neq 0$ at any point of $M^{2}$, and the level curves of $G$ are circles in $M^{2}$ with constant curvature}
$\kappa=\frac{3}{8}|\nabla G|/(c-G);$
\item [(iii)] $G$ {\it satisfies}
$$(c-G)\Delta G-|\nabla G|^2-\frac{8}{3}G(c-G)^2=0.$$
\end{itemize}
\end{theorem}

The following  uniqueness result for biconservative surfaces is due to D. Fetcu, S. Nistor and C. Oniciuc \cite{FNO16}.

\begin{theorem}\label{T:5.14} {\it Let $(M^2,g)$ be an abstract surface and $c\in\mathbb{R}$. If $M^{2}$ admits two biconservative immersions into $R^3(c)$ such that the gradients of their mean curvature functions are everywhere nonzero, then the two immersions differ by an isometry of $R^{3}(c)$.}
\end{theorem}

Even if the notion of a biconservative submanifold belongs, obviously, to extrinsic geometry, in the particular case of biconservative surfaces in $N^3(c)$ one can give an intrinsic characterization of such surfaces.

The next result is also due to D. Fetcu, S. Nistor and C. Oniciuc \cite{FNO16}.

\begin{theorem}\label{T:5.15} {\it If $\left(M^2,g\right)$ is an abstract surface, then  locally  $M^{2}$ can be isometrically embedded in  $R^3(c)$ as a biconservative surface such that $\nabla\alpha \ne 0$ everywhere on $M^{2}$ if and only if the Gauss curvature  satisfies $G<c$ and $\nabla G\neq 0$ everywhere on $M^{2}$, and its level curves are circles in $M^{2}$ with constant curvature
$\kappa=\frac{3}{8} | \nabla G |/(c-G).$}
\end{theorem}

In \cite{Ni-Thesis}, S. Nistor also  proved  the next result.

\begin{theorem} \label{T:5.16} {\it Let $\left(M^2,g\right)$ be an abstract surface and $c\in\mathbb{R}$. Assume that $G<c$ and $\nabla G\neq0$ everywhere on $M^{2}$, and the level curves of $G$ are circles in $M^{2}$ with constant curvature $\kappa=\frac{3}{8} | \nabla G |/(c-G).$
Then, locally, there is a unique biconservative embedding $\phi:\left(M^2,g\right)\to R^3(c)$ such that the mean curvature function is positive and its gradient is nowhere zero.}
\end{theorem}

\subsection{Uniqueness of complete biconservative surfaces  in $R^{3}(c)$}\label{S5.6}

 In \cite{NO22}  S. Nistor and C. Oniciuc constructed complete, simply connected, non-CMC biconservative surfaces in $R^{3}(c)$ from either the extrinsic or the
intrinsic point of view, and they proved the following uniqueness results.

\begin{theorem} \label{T:5.17} {\it  If $\phi: M\to \mathbb E^{3}$ is a complete, simply connected,  non-CMC
biconservative surface in $\mathbb E^{3}$, then, up to isometries of the domain and codomain, 
$M=\(\mathbb R^{2}, g=c^{2}(\cosh^{6} x) (dx^{2}+dy^{2})\)$ and 
$$\phi(x,y)=\(\frac{c}{3}(\cosh^{3}x)\cos(3y), \frac{c}{3}(\cosh^{3}x)\sin(3y), \frac{c}{2}\(\frac{1}{2} \sinh (2x)+x\)\),$$
where $c$ is a positive real constant.}
\end{theorem}

\begin{theorem} \label{T:5.18} {\it Let $\Phi : M^{2}\to S^{3}(1)$ be a complete,  simply connected, non-CMC biconservative surface in $S^{3}(1)$. Then, up to isometries of the domain and codomain, $M^{2}$ and
$\Phi$ are those given in Theorem 3.11 of \cite{NO22}.}
\end{theorem}

\begin{theorem} \label{T:5.19} {\it Let $\Phi : M^2 \to H^{3}(-1)$ be a complete,  simply connected, non-CMC biconservative surface in $H^{3}(-1)$. Then, up to isometries of the domain and codomain, $M^{2}$ and
$\Phi$ are those given in Theorem \ref{T:5.9}.}
\end{theorem}

\begin{theorem} \label{T:5.20} {\it Every  compact biconservative surface in $R^{3}(c)$,  $c\in \{-1,0\}$,  is CMC.}
\end{theorem}

\section{Biconservative hypersurfaces in $\mathbb E^{n+1}$}\label{S6}

\subsection{First result on biconservative hypersurfaces}\label{S6.1}

The following result of Hasanis and Valchos \cite{HV95} obtained in 1995 eas the first result on biconservative submanifolds.

\begin{theorem}\label{T:6.1} {\rm \cite{HV95}} {\it Let $M$ be a biconservative hypersurface of  $\mathbb E^4$. Then $M$ is congruent to one of the following hypersurfaces:}
\begin{enumerate}
 \item[{\rm (a)}] {\it $M$ is CMC;}

\item[{\rm (b)}] {\it $M$ is a CMC rotational hypersurface generated by a unit speed planar curve $\gamma=(f(s),g(s))$, where $f$ satisfies} $$3ff''=2(1-(f')^2);$$

\item[{\rm (c)}] {\it A generalized cylinder on a surface of revolution in $\mathbb E^3$ with non-constant mean curvature given by
$$\phi(u,v,s)=\big(f(s)\cos u,f(s)\sin u,v,g(s)\big),$$
where $\gamma=(f(s),g(s))$ is a unit speed curve with non-constant curvature and $f$ satisfies} $$3ff''=1-(f')^2;$$

\item[{\rm (d)}] {\it A $SO(2)\times SO(2)$-invariant hypersurface of $\mathbb E^{3}$ with non-constant mean curvature given by
$$\big(f(s)\cos u,f(s)\sin u,g(s)\cos v,g(s)\sin v\big),$$
where $\gamma=(f(s),g(s))$ is a unit speed curve with non-constant curvature and $f$ satisfies} $$f'g''-f''g'=\text{\small${1\over 3}\(\frac{f'}{g}-\frac{g'}{f}\)$}.$$
\end{enumerate}\end{theorem}

In \cite{MOR16},  Montaldo,  Oniciuc, and  Ratto performed a detailed qualitative study of biconservative hypersurfaces that are $SO(p+1)\times SO(q+1)$-invariant in $\mathbb E^{p+q+2}$, or $SO(p+1)$-invariant in $\mathbb E^{p+2}$. They obtained the following two results.

\begin{theorem}\label{T:6.2}  {\it There exists an infinite family of $SO(p+1)\times SO(q+1)$-invariant proper biconservative hypersurfaces (cones) in ${\mathbb E}^{p+q+2}$, and 
their corresponding profile curves $\gamma(s)$  tend asymptotically to the profile of a minimal cone. 

In addition, If $p+q \leq 17$, at infinity the profile curves $\gamma$ intersect the profile of the minimal cone at infinitely many points, while, if $p+q \geq 18$,  at infinity the profile curves $\gamma$ do not intersect the profile of the minimal cone.}
\end{theorem}

\begin{theorem} \label{T:6.3}{\it There exists an infinite family of complete $SO(p+1)$-invariant proper biconservative hypersurfaces in $\mathbb E^{p+2}$.
In addition, their corresponding profile curves  are of {catenary} type.}
\end{theorem}

\subsection{Biconservative hypersurfaces in $\mathbb E^{n+1}$}\label{S6.2}

For biconservative hypersurface  in  $\mathbb E^{5}$,    R. S. Gupta {\rm  \cite{GA19}} proved the following.

\begin{theorem}\label{T:6.5}  {\it Every biconservative hypersurface in $\mathbb E^{5}$  with constant norm of second fundamental form has constant mean curvature.}
\end{theorem}

N. C. Turgay investigated in \cite{Turgay15} biconservative hypersurfaces with three distinct curvatures in Euclidean spaces of arbitrary dimension. He proved such a hypersurface is either a generalized rotational hypersurface  or a  generalized cylinder over a rotational hypersurface. Moreover, he obtained the complete
classification of such biconservative hypersurfaces. In addition, he constructed explicit examples of such biconservative hypersurfaces in a Euclidean space.

 The next result on CMC biconservative surface is due  to  Montaldo,  Oniciuc and  Ratto {\rm  \cite{MOR16.2}}.
 
 \begin{theorem} \label{T:6.6} {\it  Let $\phi: M^{2}\to \mathbb E^{4}$ be a proper biconservative surface with constant mean curvature $\alpha\ne 0$. If $M^{2}$ has non-parallel mean curvature vector, then  locally the surface is given by
\begin{equation}\label{6.1} \phi(u, v) = (\gamma(u), v + a) = (\gamma_{1}(u), \gamma_{2}(u), \gamma_{3}(u), v + a),\;\;  a\in  \mathbb R , \end{equation}
 where $\gamma: I\to \mathbb E^{3}$ is a unit speed curve with nonzero constant 
curvature, and non-zero torsion. 
  Conversely, a surface defined by \e{6.1} is a proper biconservative surface with nonzero constant mean curvature.}
\end{theorem}

 \section{$\delta(r)$-ideal biconservative hypersurfaces in $\mathbb E^{n+1}$}\label{S7}
 
Let $M$ be a Riemannian $n$-manifold with $n\geq 3$. For a plane section $\pi\subset T_pM$, $p\in M$,  let  $K(\pi)$ denote the sectional curvature of $M$ associated with a plane section $\pi$. 
For an $s$-dimensional plane section of $L_{s}\subset T_pM$ with $s\geq 2$,  the scalar curvature $\tau(L_{s})$ of $L$ is defined by
\begin{equation}\label{7.1}\tau(L_{s})=\sum_{\alpha<\beta} K(E_\alpha\wedge E_\beta),\;\; 1\leq \alpha,\beta\leq s\end{equation}
where  $\{E_1,\ldots,E_s\}$ is an orthonormal basis of $L_{s}$. 

For any integer  $r\in [2,n-1]$, the  author defined the {\it $\delta$-invariant}, $\delta(r)$, of $M$ at $p\in M$ by  (see \cite{Chen98,Chen00,book11})
\begin{equation}\label{7.2} \delta(r)(p)=\tau(p)- \inf\{\tau(L_{r})\},\end{equation} where $L_{r}$ run over all $r$-dimensional linear subspaces of $T_pM$. 
For any $n$-dimensional submanifold $M$ of $\mathbb E^m$ and any integer $r\in [2, n-1]$,  the author established the following general optimal inequality  (cf. \cite{Chen98,Chen00,book11}):
\begin{equation}\label{7.3} \delta(r) \leq  \frac{n^2(n-r)}{2(n-r+1)} \|H\|^2.\end{equation}

\begin{definition} An $n$-dimensional  submanifold  in $\mathbb E^m$ is called {\it $\delta(r)$-ideal} if it satisfies the equality case of \eqref{7.3} identically. 
\end{definition}
Roughly speaking an ideal submanifold is a submanifold which receives the least possible tension from its ambient space. For more detailed on the theory of $\delta$-invariants, we refer to \cite{Chen93,Chen98,Chen00,book11,Chen13,CDVV13}.
 
 In 2013, the author and M. I. Munteanu proved the following.
 
 \begin{theorem}\label{T:7.1} {\rm \cite{CM13}} {\it Every $\delta(2)$-ideal biconservative hypersurface  in Euclidean spaces 
 $\mathbb E^{n+1}\,(n\geq 3)$ has constant mean curvature.}
\end{theorem}
 
In \cite{DA18}, Deepika and A. Arvanitoyeorgos improved Theorem \ref{T:7.1} to the following.

\begin{theorem}\label{T:7.2} {\rm \cite{DA18}} {\it Every $\delta(2)$-ideal  biconservative hypersurface  in  $\mathbb E^{n+1}$ is minimal.}
\end{theorem}

\begin{theorem} \label{T:7.3} {\rm \cite{DA18}} {\it Every $\delta(3)$-ideal  biconservative hypersurface  in  $\mathbb E^{5}$ has constant mean curvature.}
\end{theorem}

\begin{theorem} \label{T:7.4} {\rm \cite{DA18}} {\it Every $\delta(4)$-ideal  biconservative hypersurface in  $\mathbb E^{5}$ with constant scalar curvature has constant mean curvature.}
\end{theorem}

In 2020, Deepika and  Arvanitoyeorgos proved the following.

\begin{theorem}\label{T:7.5} {\rm \cite{DA20}} {\it Every $\delta(r)$-ideal oriented biconservative hypersurface with at most $r+1$ distinct principal curvatures in  $\mathbb E^{n+1}\,(n\geq 3)$ has constant mean curvature.}
\end{theorem}

\section{PNMCV biconservative surfaces in $\mathbb E^{4}$ and $S^{4}$}\label{S8}

A submanifold of a Riemannian manifold is said to be {\it PNMCV\/} if its mean curvature vector $H$ is nowhere zero and also its normalized mean curvature vector $\zeta = H/\|H\|$ is parallel in the normal bundle, i.e., $D\zeta=0$. 
A PNMCV  submanifold  is called {\it proper} if the mean curvature $\alpha=\| H\|$  and its gradient $\nabla \alpha$ is nowhere vanishing.

The study of submanifolds with parallel normalized mean curvature vector (PNMCV) was initiated in 1980 (see \cite{chen80}). 
In \cite{chen80}, the author proved the following results for PNMCV surfaces.

\begin{theorem}\label{T:8.1}  {\it Let $M$ be an analytic surface in a complete, simply-connected
real space form $R^{m}(c)$. If $M$ has PNMCV,  then either $M$ lies in a
hypersphere of $R^{m}(c)$ as a minimal surface or $M$ lies in a $4$-dimensional
totally geodesic submanifold $R^{4}(c)$ of $R^{m}(c)$.}
\end{theorem}

\begin{theorem}\label{T:8.2}  {\it Let $M$  be a flat analytic surface in $\mathbb E^{m}$. If $M$ has PNMCV, then $M$ is one of the following surfaces:}
\begin{enumerate}
\item[{\rm (a)}] {\it a fiat minimal surface of a hypersphere $S^{m-1}$ of $\mathbb E^{m}$;}
\item[{\rm (b)}]  {\it an open piece of the product surface of two plane circles; or}
\item[{\rm (c)}] {\it a developable surface in a linear $3$-space $\mathbb E^{3}\subset \mathbb E^{m}$.}
\end{enumerate}\end{theorem}

\begin{theorem}\label{T:8.3}  {\it Let $M$ be a compact oriented analytic surface of genus
zero in $\mathbb E^{m}$. If $M$ has PNMCV, then either M lies in a hypersphere $S^{m-1}$ of  $\mathbb E^{m}$ as a minimal
surface or $M$ lies in a linear 3-space $\mathbb E^{3}\subset \mathbb E^{m}$.}
\end{theorem}

Ye\u{g}in \c{S}en and Turgay classified in  \cite{ST18} proper PNMCV biconservative surfaces in $\mathbb E^{4}$ as follows. 

\begin{theorem} \label{T:8.4} {\it Let $M$ be a proper PNMCV biconservative surface in $\mathbb E^{4}$. Then, locally, it is congruent to the rotational surface:
$$\phi(s, t) = (\gamma_{1} \cos t, \gamma_{1} \sin t, \gamma_{2}, \gamma_{3}) $$ where  $\gamma(s)=(\gamma_{1}(s),\gamma_{2}(s),\gamma_{3}(s))$ is 
a unit speed the profile curve of $M$ such that $\gamma_{1}(s) = 1/(c \alpha^{3/4})$ with constant $c$ and $\alpha=\| H\|$, and  the curvature $\kappa$ and torsion $\tau$ of $\gamma$ are given by
 $$\kappa=\alpha(s) \sqrt{1+c^{2}\alpha(s)}, \;\; \tau(s)=\frac{c_{2}\alpha'(s)}{2\sqrt{\alpha(s)(1+c^{2}\alpha(s))}}.$$}
\end{theorem}

 In \cite{S22}, Ye\u{g}in \c{S}en studied PNMCV biconservative submanifolds $M^{n}$ of  in $\mathbb E^{n+2}$ whose shape operator has at most two distinct principal curvatures in the direction of the mean curvature vector. And she determined the canonical forms of the shape operator for such biconservative submanifolds. As an application, she shown that such biconservative submanifolds cannot be biharmonic.

Clearly, every PNMCV biconservative submanifolds of codimension 2 has flat normal bundle.
In \cite{ST18},  Ye\u{g}in \c{S}en and Turgay   investigated  meridian surfaces which lie on a rotational hypersurface of $\mathbb E^{4}$. In the same article, they constructed a family of biconservative surfaces with flat normal bundle and non-parallel normalized mean curvature vector.

PNMCV biconservative surfaces in $S^{4}$ were  studied by Nistor, Oniciuc, Turgay, and Ye\u{g}in \c{S}en in \cite{NOT23}. Among others,  they proved the following results.

\begin{theorem}  \label{T:8.5} {\rm \cite{NOT23}} {\it Any two PNMCV biconservative immersions of an abstract surface  into ${S}^4(1)$ differ by an isometry of ${S}^4(1)$.}
\end{theorem}

\begin{theorem} \label{T:8.6} \cite{NOT23}  {\it  Locally, an abstract surface $(M^2,g)$ admits a PNMCV biconservative embedding into ${S}^4(1)$ if and only if $G<1$, $\nabla G\neq 0$ everywhere, and the level curves of $G$ are circles of $M^2$ with positive constant signed curvature $\kappa=-\frac{1}{4}{\left|\nabla G\right|}/(G-1+f^2),$
where $f$ is the positive solution of the equation
$G=1-3f^2-c^2f^3$
 for some non-zero real number $c$.}
\end{theorem}

\begin{theorem}  \label{T:8.7} {\rm \cite{NOT23}} {\it If a  Riemannian $3$-manifold  admits two proper biconservative PNMCV
isometric immersions into $\mathbb E^{4}$ with three distinct principal curvatures, then the two immersions differ by an isometry of $\mathbb E^{4}$.}
\end{theorem}

Recently,  \c{S}en and Turgay \cite{ST23} investigated  PNMCV biconservative submanifolds $M^{3}$ in $\mathbb E^{5}$, and  proved that the principal curvatures and principal directions of $M^{3}$ can be determined intrinsically.
In addition, \c{S}en and Turgay \cite{ST23} obtained the following.

\begin{theorem}  \label{T:8.9}  {\it There do not exist proper biharmonic submanifolds in $\mathbb E^{5}$ with parallel normalized mean curvature vector.}
\end{theorem}

\section{Biconservative submanifolds in $R^{m}(c)$}\label{S9}

\subsection{Biconservative hypersurfaces in $R^{n+1}(c)$}\label{S9.1}

In \cite{MOP23}, Montaldo, Oniciuc and  P\'ampano  characterized the profile curves of non-CMC biconservative rotational hypersurfaces in a real space form $R^{n+1}(c)$, $c\in \{-1,0,1\},$ as $p$-elastic curves, for a rational number $p \in [1/4, 1)$ depending on $n$. They proved the existence of a discrete biparametric family of non-CMC closed biconservative hypersurfaces in $S^{n+1}(c)$.  In addition, they obtained the following two results.

\begin{theorem}\label{T:9.1} {\rm \cite{MOP23}} {\it Let $M^{n}$ be a non-CMC closed biconservative rotational hypersurface of a space form $R^{n+1}(c)$. Then $R^{n+1}(c)$  is the $n$-dimensional round sphere $S^{n+1}(c)$.}
\end{theorem}

For compact biconservative hypersurfaces in a space form $R^{m+1}(c)$, D. Fetcu, E. Loubeau and C.  Oniciuc  proved the next result in \cite{FLO21}.

\begin{theorem}\label{T:9.2} {\it Let $\phi:M\to R^{n+1}(c)$ be a compact biconservative hypersurface of a real space form $R^{n+1}(c)$, with $c\in\{-1,0,1\}$. If $M$ is non-minimal and non-negatively curved, and has constant scalar curvature, then $\phi(M)$ is either}
\begin{enumerate}

\item[(1)] {\it ${S}^n(r)$, $r>0$, whenever $c\in\{-1,0\}$; or}

\item[(2)] {\it ${S}^n(r)$, $0<r<1$, or ${S}^{n_1}(a)\times {S}^{n-n_1}(b)\subset {S}^{n+1}(1)$ with $a^2+b^2=1$ and $a\neq\sqrt{n_1/n}$, whenever $c=1$.}
\end{enumerate}
\end{theorem}

R. S. Gupta and A. Arvanitoyeorgos  studied in \cite{GA22} biconservative hypersurfaces in
real space forms with four distinct principal curvatures whose second
fundamental form has constant norm. They obtained the following.

\begin{theorem} \label{T:9.3}  {\it Every biconservative hypersurface  in $R^{n+1}(c)$
with at most four distinct principal curvatures, whose second fundamental form has constant norm, is of constant mean curvature and of constant scalar curvature.}
\end{theorem}

For compact non-minimal biconservative hypersurfaces in $R^{n+1}(c)$, we have the following result from {\rm \cite{FO22}}. 

\begin{theorem}\label{T:9.4} {\it Let $\phi:M\to  R^{n+1}(c)$ be a compact non-minimal biconservative hypersurface. If $M$ is non-negatively curved and $n\leq 10$, then $\phi$ is CMC and $\phi(M)$ is one of the hypersurfaces given by Theorem \ref{T:9.2}.}
\end{theorem}

\subsection{Biconservative hypersurfaces in $S^{4}$}\label{S9.2}

N. C. Turgay and A. Upadhyay  \cite{TU19} studied biconservative hypersurfaces in $R^{n+1}(c)$ and obtained their parametric equations. Moreover, they classified biconservative hypersurfaces  with three distinct principal curvatures in $R^4(c)$  explicitly for $c\neq 0$ as follows.
 
 \begin{theorem} \label{T:9.5} {\rm \cite{TU19}} {\it If $M$ is a proper biconservative hypersurface in 
 $S^4(1)\subset \mathbb E^{5}$ with three distinct principal curvatures at every point, then, locally, $M$ is congruent to the hypersurface in $\mathbb E^{5}$ given by
\begin{equation}\label{9.1}  \phi(s,t,u)=\big(\gamma_1 (s), \gamma_2 (s) \cos t,\gamma_2 (s)\sin t, \gamma_3(s) \cos u, \gamma_3(s)  \sin u\big)\end{equation}
 for a unit speed curve curve $\gamma=(\gamma_1 ,\gamma_2 ,\gamma_3)$ in $S^2(1)\subset \mathbb E^{3}$.}
\end{theorem}

\begin{theorem}\label{T:9.6} {\rm \cite{TU19}} {\it If $M$ is a hypersurface in $S^{4}(1)$ defined by \e{9.1} for a unit speed  curve $\gamma:(a,b)\to S^{2}(1)\subset \mathbb E^{3}$,  then $M$ is biconservative if and only if, locally, it is congruent to the hypersurface given by \e{9.1}  for the curve $\gamma$ described by 
$$\gamma(s)=(\cos u(s), \sin u(s)\cos v(s), \sin u(s)\sin v(s))$$  satisfying the following two differential equations:
\begin{equation*} v'^{2} \sin^{2} u+u'^{2}=1\end{equation*} and}
\begin{equation*} 3 u''\sin u+2 u' \cos 2v\sqrt{1-u'^{2}}=5\cos u(1-u'^{2}).\end{equation*}
\end{theorem}

 Montaldo,  Oniciuc and  P\'ampano \cite{MOP23} proved the following.

\begin{theorem}\label{T:9.7}  {\it For an integer $n\geq 2$, there exists a discrete biparametric family of closed non-CMC biconservative rotational hypersurfaces $M^{n}$ in ${S}^{n+1}(1)$. However, none of these hypersurfaces are embedded in ${S}^{n+1}(1)$.}
\end{theorem}

\subsection{Biconservative hypersurfaces in $H^{4}$}\label{S9.3}

The next result is due to N. C. Turgay and A. Upadhyay.

\begin{theorem} \label{T:9.8} {\rm \cite{TU19}} {\it If $M$ is a proper biconservative hypersurface in $H^{4}(-1)$ with three distinct principal curvatures, then, locally, 
$M$ is  congruent to one of the following hypersurfaces:}
\begin{enumerate}
\item[{\rm (1)}] {\it A hypersurface in $\mathbb H^4(-1)$ defined by
$$\phi(s,t,u)=\(\gamma_{1}(s),\gamma_{2}(s)\cos t,\gamma_{2}(s)\sin t,\gamma_{3}(s)\cos u, \gamma_{3}(s)\sin u\)$$
for a unit speed curve $\gamma=(\gamma_{1},\gamma_{2},\gamma_{3}):(a,b)\to H^{2}(-1)\subset \mathbb E^{3}_{1}$;}

\item[{\rm (2)}] 
{\it A hypersurface in $\mathbb H^4(-1)$ defined by
\begin{equation*}
\phi(s,t,u)=\left(\gamma_1(s)\cosh t,\gamma_1(s)\sinh t,\gamma_2(s)\cos t,\gamma_2(s)\sin t,\gamma_3(s)\right)
\end{equation*} 
for an unit speed curve $\gamma=(\gamma_1,\gamma_2,\gamma_3):s\in (a,b)\to  H^2(-1)\subset \mathbb E^{3}_{1}$;}

\item[{\rm (3)}] 
{\it A hypersurface in $\mathbb H^4(-1)$ defined by
\begin{equation}\begin{aligned}\notag
x(s,t,u)=&\left(\frac{a \beta(s)^2+a}{s}+a s u^2+\frac{s}{4 a},s u,\beta(s) \cos t,\beta(s) \sin t,\right.\\&\hskip.3in \left.
\frac{a \beta(s)^2+a}{s}+a s u^2-\frac{s}{4 a}\right)
\end{aligned} \end{equation} 
for a function $\beta$ and some non-zero constants $a$;}

\item[{\rm (4)}] 
{\it A hypersurface in $\mathbb H^4(-1)$ defined by
\begin{equation}\begin{aligned}\notag
x(s,t,u)=&\left(\frac{a \beta(s)^2}{s}+a s \left(t^2+u^2\right)+\frac{s}{4 a}+\frac{a}{s},s t,s u,\beta(s),\right.\\&\hskip.3in  \left.
\frac{a \beta(s)^2}{s}+a s \left(t^2+u^2\right)-\frac{s}{4 a}+\frac{a}{s}\right)
\end{aligned} \end{equation}
for a function $\beta$ and a non-zero constant $a$.}
\end{enumerate}
\end{theorem}

\subsection{Biconservative surfaces with codimension $\geq 2$}\label{S9.4}

For CMC biconservative surfaces of real space form with  codimension $\geq 2$,
Montaldo,  Oniciuc and Ratto obtained the following.

\begin{theorem}\label{T:9.9} {\rm \cite{MOR16.2}} {\it If $M^{2}$ is a CMC biconservative surface in a  real space form $R^{m}(c)$ with $c \ne 0$, then $M^{2}$  has PMCV.}
\end{theorem}

\subsection{Biconservative hypersurfaces with constant scalar curvature}\label{S9.5}

In \cite{FO22}, D. Fetcu and C. Oniciuc proposed the following open problem:

\vskip.1in
\noindent {\bf Open Problem:} Classify all biconservative hypersurfaces with constant
scalar curvature in  $R^{n+1}(c)$.

\vskip.1in
 In \cite{FHZ21.2}, Y. Fu, M.-C.  Hong, and X. Zhan obtained the following two results:

\begin{theorem}\label{T:9.10} {\it Any biconservative hypersurface with constant scalar curvature in $R^4(c)$ is CMC.}
\end{theorem}

\begin{theorem}\label{T:9.11} {\it Any biconservative hypersurface with constant scalar curvature in
$R^5(c)$ is ether an open part of a certain rotational hypersurface
or a CMC hypersurface.} 
\end{theorem}

These two results solved Fetcu and C. Oniciuc's  open problem  for $n\leq4$.

\section{Biconservative Submanifolds in $S^{n}\times \mathbb R$ and $H^{n}\times \mathbb R$}\label{S10}

 Fetcu, Oniciuc, and  Pinheiro  \cite{FOP15} studied non-minimal biconservative surfaces with parallel mean curvature vector field $H$ in the product manifolds ${S}^n\times\mathbb{R}$ and ${H}^n\times\mathbb{R}$. 
 In particular,  they obtained the following classification result.

\begin{theorem}\label{T:10.1} {\rm \cite{FOP15}} {\it If $M$ is a PMCV biconservative surface with  $H\ne 0$ in $R^n(c)\times\mathbb{R}$, $c\in\{-1, 1\}$, then it is one of the following surfaces:}
\begin{enumerate}
\item[{\rm (1)}] {\it a minimal surface of a totally umbilical hypersurface of $R^n(c)$; or }

\item[{\rm (2)}] {\it a CMC surface in a $3$-dimensional totally umbilical submanifold of $R^n(c);$ or}

\item[{\rm (3)}] {\it  a vertical cylinder over a circle in $R^2(c)$ with curvature $\kappa=2\|H\|;$ or}

\item[{\rm (4)}] {\it  a surface lies in ${S}^4\times\mathbb{R}\subset\mathbb{E}^5\times\mathbb{R}$ and, as a surface in $\mathbb{E}^5\times\mathbb{R}$,  it is  locally given by 
\begin{equation}\begin{aligned}\notag \phi(u,v)=\,&\frac{1}{a}\{c_1+\sin\theta(v_1\cos(au)+v_2\sin(au))\}+(u\cos\theta+b)\xi  
\\& +\frac{1}{\kappa}(c_2(\cos v-1)+c_3\sin v),\;\; \theta\in(0,\pi/2), \end{aligned} \end{equation}
where $a=\sqrt{1+\sin^2\theta}$, $b\in \mathbb R$, $\kappa=\sqrt{1+4\|H\|^2+\sin^2\theta}$,  $\{c_2,c_3\}$ are orthonormal constant vectors in $\mathbb{E}^5\times\mathbb{R}$ such that $c_2\perp\xi$ and $c_3\perp\xi$, $c_1$ is a unit constant vector satisfying $\langle c_1,c_2\rangle=a/\kappa\in(0,1)$, $c_1\perp c_3$, $c_1\perp\xi$, and $\{v_1,v_2\}$ are   orthonormal constant vectors lying in the orthogonal complement of ${\rm Span}
\{c_1,c_2,c_3,\xi\}$ in $\mathbb{E}^5\times\mathbb{R}$.}
\end{enumerate}
\end{theorem}

For a given  isometric immersion $\phi:M\to R^{n}(c)\times \mathbb R$, let $\partial_{t}$ denote a unit  vector field tangent to the second factor. Let us define a vector  field $T$ tangent  to $M$ and a normal vector field $\eta$ of  along $\phi$ by 
\begin{equation}\label{10.1}
\partial_{t}=\phi_\ast T+\eta.
\end{equation}
Now, consider  an oriented minimal surface $$\psi:M^2\to R^{2}(a)\times \mathbb R$$
such that the vector field $T$  is  nowhere zero, where $a\neq 0$ and $|a|<1$. 
 Let $b$ be a positive number satisfying $a^2+b^2=1$. 
Consider
\begin{equation}\label{10.2}\phi:M^3=M^2\times I\to R^{4}(c)\times \mathbb R, \;\; c=\pm 1,\end{equation}
defined by
\begin{eqnarray}\label{eq:localf}\notag
\phi(p,s)=\left(b\cos\frac{s}{b},b\sin\frac{s}{b},\psi(p)\right).
\end{eqnarray}

For the map $\phi$ defined by \e{10.2},   Manfio,  Turgay and  Upadhyay proved the following.

\begin{theorem}\label{theo:main} {\rm \cite{MTU19}}
{\it The map $\phi$ defines, at regular points, an isometric immersion satisfying 
$\<H,\eta\>=0$, where $H$ is the mean curvature vector  of $\phi$.
Moreover, $\phi$ is a biconservative isometric immersion
with PMCV if and only if $\psi$ is a vertical cylinder. 

Conversely, every biconservative submanifold $\phi:M^3\to R^{4}(c)\times \mathbb R$ with nonzero PMCV such that  $T$ given by \eqref{10.1} is nowhere 
zero is locally obtained in this way.}
\end{theorem}

\section{Biconservative submanifolds  in Riemannian manifolds}\label{S11}

\subsection{General results for biconservative submanifolds}\label{S11.1}

Let $\phi:M^n\to \tilde M^{m}$ be a submanifold. Computing $\tau_2(\varphi)$ by splitting it in the tangent and in the normal part and using \eqref{subvarietateoarecare1}  one obtains the following characterizations for biconservative submanifolds (see various expressions for $\tau_2(\varphi)$ in \cite{chen83,LMO08,Ni17,Ou10,O02}).

\begin{proposition}\label{P:11.1} {\it Let $\phi:M^n\to \tilde M^{m}$ be a submanifold. Then the following conditions are equivalent:}
\begin{enumerate}
    \item[{\rm (a)}] {\it $M$ is biconservative;}
    \item[{\rm (b)}] $\trace A_{\nabla^\perp_{\cdot} H}(\cdot)+\trace (\nabla A_H) +\trace ({\tilde R}(\,\cdot\, ,H)\,\cdot\,)^T=0$;
    \item[{\rm (c)}] $\frac{n}{2}\nabla \|H\|^2+2\,\trace (A_{\nabla^\perp_{\cdot} H}(\cdot)) + 2\, \trace ({\tilde R}(\,\cdot\, ,H)\,\cdot\,)^T=0$;
    \item[{\rm (d)}] $2\,\trace (\nabla A_H)-\dfrac{n}{2}\nabla \|H\|^2=0$,
\end{enumerate}
{\it where ${\tilde R}$ denotes the curvature tensor of $\tilde M^{m}$.}
\end{proposition}

As an easy consequence of Proposition \ref{P:11.1}, we have

\begin{proposition}\label{P:11.2} {\it Let $M^{n}$ be a PMC submanifold of a Riemannian manifold $\tilde M^{m}$. Then $M^{n}$ is biconservative if and only if }
$$\trace\,(\overline{R}(\,\cdot\,,H)\,\cdot\,)^{T}=0.$$
\end{proposition}

\begin{proposition}\label{P:11.3} {\rm \cite{Ni17}} {\it Let $\phi:M^n\to \tilde M^{m}$ be a submanifold. Then we have:}
\begin{enumerate}
  \item[{\rm (1)}]  {\it the stress-bienergy tensor of $\varphi$ is determined by}
   \begin{equation}\label{11.1}
    S_2= - \frac{n^2}{2}\|H\|^2 I + 2n A_H;
   \end{equation}
  \item[{\rm (2)}] $\trace \,S_2 = \left(2-\dfrac{n}{2}\right)n^2 \|H\|^2$;
  \item[{\rm (3)}] {\it the relation between the divergence of $S_2$ and the divergence of $A_H$ is given by}
  \begin{equation}\label{11.2}
  \Div \, (S_2)=-\frac{n^2}{2}\nabla \left(|H|^2\right)+2n\,\Div (A_H);
  \end{equation}
  \item[{\rm (4)}] $\|S_2\|^2=\left(\dfrac{n}{4}-2\right)n^4 \|H\|^4+4n^2\|A_H\|^2$.
\end{enumerate}
\end{proposition}

\begin{remark} From  \e{11.2} we see that ``$M$ is biconservative''  does not imply ``$\Div (A_H)=0$''. In fact, only under the condition $\|H\|$ is constant, the biconservativity is equivalent to $\Div (A_H)=0$.
\end{remark}

\subsection{Biconservative surfaces in Riemannian manifolds}\label{S11.2}

For an isometric immersion $\phi:M^{n} \to \tilde M^{m}$, let $\phi_{H}$  denote the  traceless part of $A_H$, i.e.,
\begin{equation} \phi_H= A_{H}-\| H\|^{2}I,\end{equation} 
where $I$ is the identity map.

For non-minimal CMC biconservative surfaces,  Fetcu, Oniciuc and    Pinheiro \cite{FOP15} proved the following.

\begin{theorem}\label{T:11.5}  \cite{FOP15}
{\it If $M$ is a non-minimal CMC biconservative surface in a Riemannian manifold $\tilde M$, 
then the Gauss curvature  of $M$ satisfies}
$$\frac{1}{2}\Delta\|\phi_H\|^2=2\|\phi_H\|^2G +\|\nabla\phi_H\|^2,$$
\end{theorem}

Theorem \ref{T:11.5} implies the following.

\begin{corollary}\label{C:11.6}   \cite{FOP15} {\it Let $M$ be a CMC biconservative surface in a Riemannian manifold $\tilde M$. If $M$ is compact and $G\geq 0$, then $\nabla A_H=0$ and $M$ is either pseudo-umbilical or flat.}
\end{corollary}

By definition, the \textit{principal curvatures} of a submanifold $M$ of  Riemannian manifold $\tilde M$ are defined to be the eigenvalue functions of $A_H$.

The next three results were due to S. Nistor.

\begin{proposition} \label{P:11.7}  {\rm \cite{Ni17}} {\it Let $\phi: M\to \tilde M$ be a $CMC$ biconservative surface. If $M$ is compact and does not contain pseudo-umbilical points, then $M$ is a topologic torus.}
\end{proposition}

\begin{proposition}\label{P:11.8} {\rm \cite{Ni17}}  {\it Let $\phi: M\to \tilde M$ be a $CMC$ biconservative surface such that $M$ is compact and does not contain pseudo-umbilical points. If $G\geq 0$ or $G\leq 0$, then $\nabla A_H=0$ and $M$ is flat}
\end{proposition}

\begin{theorem} \label{T:11.9}  {\rm \cite{Ni17}}  {\it Let $\phi: M\to \tilde M$ be a biconservative surface. If the principal curvatures  $\lambda_{1},\lambda_{2}$ of $M$ are constant, the $\nabla A_{H}=0$.}\end{theorem}

The following two results of {\rm \cite{Ni17}} say that $CMC$ biconservative surfaces in a Riemannian manifold  can be immersed in $R^3(c)$ which has  the shape operator either the tensor field $A_H$ or the $S_2$ as given by \e{11.1}.

\begin{theorem} \label{T:11.10}  {\rm \cite{Ni17}}  {\it Let $\phi: M\to \tilde M$ be a $CMC$ biconservative surface in a Riemannian manifold and let $\lambda_1$ and $\lambda_2$ denote the principal curvatures of $M$ corresponding to $\phi$. If $\lambda_1$ and $\lambda_2$ are constants and $\lambda_1>\lambda_2$, then we have:}
\begin{itemize}
  \item [{\rm (a)}] {\it locally, there exists $\psi:M\to R^3(c)$ an isoparametric surface such that $A^\phi_{H^\phi}$ is the shape operator of  $\psi$ (with respect to the unit normal vector field), where
  $c = {\mu^2}/{4}-\|H^\phi\|^4$
and $\mu=\lambda_{1}-
      \lambda_{2}$.  Moreover, we have $\|H^\psi \|=\|H^\phi\|^2$.}
      
  \item [{\rm (b)}] {\it locally, there exists an isoparametric surface $\psi:M\to R^3(c)$  such that $S^\phi_2$ is the shape operator of $\psi$ (with respect to the unit normal vector field), where
      $c=4\left(\mu^2-\|H^\phi\|^4\right). $
      Moreover, we have $\|H^\psi\|=2\|H^\phi\|^2$.}
\end{itemize}
\end{theorem}

\section{Biconservative hypersurfaces in BCV-spaces}\label{S12}

Three-dimensional simply-connected homogeneous manifolds have been classified according to the dimension $d$ of their isometry group which is equal to 3, 4 or 6. If $d=6$, one obtains the real space forms. 

The~Bianchi--Cartan--Vranceanu spaces (simply called {\it BCV-spaces\/}) are three-dimensional homogeneous Riemannian manifolds with isometry group of dimension 4 or 6.
Such spaces can be  described by the following 2-parameter family of Riemannian metrics,  {\it Bianchi-Cartan-Vranceanu metrics}, which are defined by
\begin{equation}\label{12.1}
g_{\kappa,\tau} =\frac{dx^{2} + dy^{2}}{F^{2}} +  \left(dz +
\tau\, \frac{ydx - xdy}{F}\right)^{2},\;\; \kappa, \tau\in\mathbb R,
\end{equation}
defined on  ${\mathcal M}=\{(x,y,z)\in{\mathbb R}^3\colon F(x,y)> 0\}$ with $F=1+\dfrac{\kappa}{4}(x^2+y^2).$

It is well-known that 
BCV-spaces  admit a Riemannian submersion, called {\it Hopf fibration}, over a surface with constant Gauss curvature $\kappa$:
$$\pi:{\mathcal M}_{\kappa,\tau}= ({\mathcal M}, g_{\kappa,\tau})\to M^2(\kappa)=\({\mathbb R}^2, g_{F}=\frac{dx^2+dy^2}{F^2}\)$$  such that $\pi(x,y,z)=(x,y)$.
The vector field $\partial _{z}=\partial/\partial z$ is tangent to the fibers of the Hopf fibration, which is known as the {\it Hopf vector field}. 

B. Daniel \cite{Da} considered  the angle $\theta$ between the normal vector field of a surface in ${\mathcal M}_{\kappa,\tau}$ and $\partial _{z}$ in terms of $\cos\theta$. Moreover, he proved that the angle function $\theta$ is a fundamental invariant for surfaces in BCV-spaces.

The 2-parameter family  of metrics given by \e{12.1} includes all the $3$-dimensional homogeneous metrics whose isometry group has dimension $4$ or $6$, except for the hyperbolic space, according to 
\begin{itemize}
\item[{\rm (a)}] ${\mathcal M}_{\kappa,\tau}\cong \mathbb E^3$,  if $\kappa=\tau=0$;
\item[{\rm (b)}]  ${\mathcal M}_{\kappa,\tau}\cong S^{3}\left(\frac{\kappa}{4}\right)\setminus\{\infty\}$,  if $\kappa=4\tau^2\neq 0$;
\item[{\rm (c)}]   ${\mathcal M}_{\kappa,\tau}\cong (S^2(\kappa)\setminus\{\infty\})\times\mathbb R$, if $\kappa >0$ and $\tau=0$;
\item[{\rm (d)}]  ${\mathcal M}_{\kappa,\tau}\cong H^2(\kappa)\times\mathbb R$, if $\kappa <0$ and $\tau=0$;
\item[{\rm (e)}]   ${\mathcal M}_{\kappa,\tau}\cong\mathrm{SU}(2)\setminus\{\infty\}$, if $\kappa >0$ and $\tau\neq 0$;
\item[{\rm (f)}]   ${\mathcal M}_{\kappa,\tau}\cong\widetilde{\mathrm{SL}(2, R)}$, if $\kappa <0$ and $\tau\neq 0$;
\item[{\rm (g)}]   ${\mathcal M}_{\kappa,\tau}\cong {\rm Nil}_3$, if $\kappa =0$ and $\tau\neq 0$.
\end{itemize}

  Montaldo, Onnis and  Passamani  \cite{MOP17} studied biconservative surfaces in 
a Bianchi-Cartan-Vranceanu spaces $\Nkt$ and proved the following results.

\begin{proposition}\label{P:12.1} {\it
	Let $M$ be a constant angle surface in $\Nkt$ with angle $\theta \in [0,\pi]$. If  $ M$ is biconservative, then it is CMC.}    
\end{proposition}

\begin{proposition}\label{P:12.2} {\it
	Let $M$ be a  non-minimal biconservative surface in $\Nkt$ with $\kappa\neq 4\tau^2$ . If $M$ is  CMC, then it is a Hopf tube over a curve with constant geodesic curvature $($via the Hopf fibration $\pi:{\mathcal M}_{\kappa,\tau}\to M^2(\kappa))$.} 
\end{proposition}

\begin{theorem}\label{T:12.3} {\it
If $M$ is a non-minimal biconservative surface in  $\Nkt$ with $\kappa\neq 4\tau^2$, then  the following three statements are equivalent:}
	\begin{itemize}
		\item[(a)] {\it $M$ is a constant angle surface;}
		\item[(b)] {\it $M$ is a CMC surface;}
		\item[(c)] {\it $ M$ is a Hopf tube over a curve with constant geodesic curvature.}
	\end{itemize}	
\end{theorem}

\begin{theorem}\label{T:12.4} {\it
Let $M$ be a surface of revolution in  $\Nkt$ which is not a real space form and with $\tau\neq0$. If $\theta\in(0,\pi) $ and the mean curvature $\alpha$ of $M$ is nowhere zero, then $M$ is a  biconservative surface if and only if $M$ is a Hopf circular cylinder. } 
\end{theorem}

\section{Biconservative submanifolds in Hadamard manifold}\label{S13}

A {\it Hadamard manifold} is a complete, simply-connected Riemannian manifold with sectional curvatures $Riem^{M}\leq 0$. 

Fetcu, Oniciuc, and  Pinheiro studied non-minimal CMC biconservative surfaces in Hadamard manifolds. They proved the following results.

\begin{theorem}\label{theo-compact0} {\rm \cite{FOP15}} {\it
Let $M$ be a complete, non-minimal, CMC biconservative surface in a Hadamard manifold $\tilde M$, with sectional curvature bounded from below by a negative constant $c$, such that the norm of its second fundamental form $M$ is bounded, i.e.,
$$\int_{M}\| \phi_H \|^2\ dv < +\infty ,$$
and $\|H\|^2>\frac{1}{2}(\zeta-2c)$ where  $\zeta=\sup_{M}
\(|\sigma|^2-\tfrac{\|A_{H}\|^2}{\|H\|^2}\).$ Then $M$ is compact.}
\end{theorem}

When $n=2$, they obtained the following.

\begin{corollary} {\rm \cite{FOP15}} {\it
If $M$ is a complete, non-minimal, CMC biconservative surface in a  Hadamard 3-manifold $\tilde M$ with sectional curvature bounded from below by a negative constant $c$ such that 
$\|H\|^2>-c$ and$$\int_{M}\|\phi_H\|^2\ dv < +\infty,$$
then $M$ is compact.}
\end{corollary}

\section{Biconservative surfaces in Lorentzian 3-space forms $R^{3}_{1}(c)$}\label{S14}

Nondegenerate biconservative surfaces in Lorentzian 3-space forms $R^{3}_{1}(c)$ were completely classified by Yu Fu in \cite{Fu15}.

\subsection{Lorentzian biconservative surfaces in de Sitter 3-space $S^{3}_{1}$}\label{S14.1}

For nondegenerate biconservative surface in the de Sitter space $S_{1}^{3}(1)$, Fu proved the following classification theorem. 

\begin{theorem}\label{S14.1} {\it If $\phi : M\to S_{1}^{3}(1)\subset \mathbb E_{1}^{4}$ is a nondegenerate biconservative surface in $S_{1}^{3}(1)$, then $M$ is either a CMC surface or locally given by one of the following nine surfaces:}
\begin{enumerate}
 \item[{\rm (1)}] {\it a spacelike rotational surface defined by
 $$ \phi(u,v)= \(u \cosh v, u\sinh v,\sqrt{1+u^{2}} \cos f ,\sqrt{1+u^{2}} \sin f\),$$
 where $u\in (0,\infty)$ and }
 $$f=\pm \int \frac{3du}{u^{1/3}(1+u^{2})\sqrt{9 u^{-2/3}-u^{2}-1}}.$$
 
 \item[{\rm (2)}] {\it a spacelike rotational surface defined by
 $$ \phi(u,v)= \(\sqrt{1-u^{2}}\sinh f, \sqrt{1-u^{2}}\cosh f, u \sin v, u\cos v\),$$
 where $u\in (0,1)$ and }
 $$f=\pm \int \frac{3du}{u^{1/3}(1-u^{2})\sqrt{1+9 u^{-2/3}-u^{2}}}.$$
 
\item[{\rm (3)}] {\it a spacelike rotational surface defined by
 $$ \phi(u,v)= \(\sqrt{u^{2}-1}\cosh f, \sqrt{u^{2}-1}\sinh f, u \sin v, u\cos v\),$$
 where $u\in (1,\infty)$ and }
 $$f=\pm \int \frac{3du}{u^{1/3}(u^{2}-1)\sqrt{1+9 u^{-2/3}-u^{2}}}.$$

\item[{\rm (4)}]  {\it a spacelike rotational surface defined by
 $$ \phi(u,v)= \(\frac{1}{2}\! \left\{u(v^{2}+ f^{2})-\frac{1}{u}+u\right\}\! ,\frac{1}{2}\!\left\{u(v^{2}+ f^{2})-\frac{1}{u}-u\right\}\! ,uv,vf \),$$
 where $u\in (0,3^{3/4})$ and }
 $$f= \int \frac{3du}{u^{2}\sqrt{9- u^{8/3}}}.$$

\item[{\rm (5)}] {\it a timelike rotational surface defined by
 $$\phi(u,v)= \(u \sinh v, \sqrt{1-u^{2}}\cos f,  \sqrt{1-u^{2}}\sin f,u \cosh v\),$$
 where $u\in (0,1)$ and}
$$f=\pm \int \frac{3du}{u^{1/3}(1-u^{2})\sqrt{1-9 u^{-2/3}-u^{2}}}.$$

\item[{\rm (6)}] {\it a timelike rotational surface defined by
 $$ \phi(u,v)= \(u \cosh v, u\sinh v,\sqrt{1+u^{2}} \cos f ,\sqrt{1+u^{2}} \sin f\),$$
 where $u\in (0,\infty)$ and }
 $$f=\pm \int \frac{3du}{u^{1/3}(1+u^{2})\sqrt{9 u^{-2/3}+u^{2}+1}}.$$

\item[{\rm (7)}] {\it a timelike rotational surface defined by
 $$ \phi(u,v)= \(\sqrt{1-u^{2}}\sinh f, \sqrt{1-u^{2}}\cosh f, u \sin v, u\cos v\),$$
 where $u\in (0,1)$ and }
 $$f=\pm \int \frac{3du}{u^{1/3}(1-u^{2})\sqrt{9 u^{-2/3}+u^{2}-1}}.$$

\item[{\rm (8)}] {\it a timelike rotational surface defined by
 $$ \phi(u,v)= \(\sqrt{u^{2}-1}\sinh f, \sqrt{u^{2}-1}\cosh f, u \sin v, u\cos v\),$$
 where $u\in (1,\infty)$ and }
 $$f=\pm \int \frac{3du}{u^{1/3}(u^{2}-1)\sqrt{9 u^{-2/3}+u^{2}-1}}.$$

\item[{\rm (9)}]   {\it a timelike rotational surface defined by
 $$ \phi(u,v)= \(\frac{1}{2}\! \left\{u(v^{2}+ f^{2})-\frac{1}{u}+u\right\}\! ,\frac{1}{2}\!\left\{u(v^{2}+ f^{2})-\frac{1}{u}-u\right\}\! ,uv,vf \),$$
 where $u\in (0,\infty)$ and }
 $$f= \int \frac{3du}{u^{2}\sqrt{9+u^{8/3}}}.$$
\end{enumerate}\end{theorem} 

\begin{remark} The surfaces given by (2), (3), (7) and (8) in Theorem \ref{S14.1} are spherical; the surfaces given by (1), (5) and (6) are hyperbolical, and  the surfaces given by (5) and (9) are parabolic.
\end{remark}

\subsection{Lorentzian biconservative surfaces in anti de Sitter 3-space $S^{3}_{1}$}\label{S14.2}

In \cite{Fu15}, Y. Fu also  obtained  the explicit classifications of nondegenerate biconservative surfaces in  $H_{1}^{3}(-1)$ as follows.

\begin{theorem}\label{S14.3} {\it If $\phi : M\to H_{1}^{3}(-1)\subset \mathbb E_{2}^{4}$ is a nondegenerate biconservative surface $H_{1}^{3}(-1)$, then $M$ is either a CMC surface or locally given by one of the following eleven surfaces:}
\begin{enumerate}
\item[{\rm (1)}] {\it a spacelike rotational surface defined by
 $$ \phi(u,v)= \(u \cosh v, u\sinh v,\sqrt{1-u^{2}} \cosh f ,\sqrt{1-u^{2}} \sinh f\),$$
 where $u\in (0,1)$ and }
 $$f=\pm \int \frac{3du}{u^{1/3}(1-u^{2})\sqrt{9 u^{-2/3}+u^{2}-1}}.$$
 
 \item[{\rm (2)}] {\it a spacelike rotational surface defined by
 $$ \phi(u,v)= \(u \cosh v, u\sinh v,\sqrt{u^{2}-1} \cosh f ,\sqrt{u^{2}-1} \sinh f\),$$
 where $u\in (1,\infty)$ and }
 $$f=\pm \int \frac{3du}{u^{1/3}(u^{2}-1)\sqrt{9 u^{-2/3}+u^{2}-1}}.$$ 
 
 \item[{\rm (3)}] {\it a spacelike rotational surface defined by
 $$ \phi(u,v)= \(\sqrt{1+u^{2}}\cos f, \sqrt{1+u^{2}}\sin f, u \sin v, u\cos v\),$$
 where $u\in (0,\infty)$ and }
 $$f=\pm \int \frac{3du}{u^{1/3}(1+u^{2})\sqrt{1+9 u^{-2/3}-u^{2}}}.$$

\item[{\rm (4)}]  {\it a spacelike rotational surface defined by
 $$ \phi(u,v)= \(\frac{1}{2}\! \left\{u(v^{2}- f^{2})+\frac{1}{u}+u\right\}\! ,uf,\frac{1}{2}\!\left\{u(v^{2}- f^{2})+\frac{1}{u}-u\right\}\! ,uv \),$$
 where $u\in (0,\infty)$ and }
 $$f= \int \frac{3du}{u^{2}\sqrt{9+u^{8/3}}}.$$

 \item[{\rm (5)}] {\it a timelike rotational surface defined by
 $$\phi(u,v)= \(u \sinh v, \sqrt{1+u^{2}}\cosh f,  \sqrt{1+u^{2}}\sinh f,u \cosh v\),$$
 where $u\in (0,\infty)$ and}
$$f=\pm \int \frac{3du}{u^{1/3}(1+u^{2})\sqrt{1-9 u^{-2/3}+u^{2}}}.$$

\item[{\rm (6)}] {\it a timelike rotational surface defined by
 $$ \phi(u,v)= \(u \cos v, u\sin v,\sqrt{u^{2}-1} \cos f ,\sqrt{u^{2}-1} \sin f\),$$
 where $u\in (1,\infty)$ and }
 $$f=\pm \int \frac{3du}{u^{1/3}(u^{2}-1)\sqrt{u^{2}-9 u^{-2/3}-1}}.$$ 

\item[{\rm (7)}]  {\it a timelike rotational surface defined by
 $$ \phi(u,v)= \(\frac{1}{2}\! \left\{u(v^{2}+ f^{2})+\frac{1}{u}+u\right\}\! ,uv,\frac{1}{2}\!\left\{u(v^{2}+ f^{2})+\frac{1}{u}-u\right\}\! ,uf \),$$
 where $u\in (3^{3/4},\infty)$ and }
 $$f= \int \frac{3du}{u^{2}\sqrt{u^{8/3}-9}}.$$

\item[{\rm (8)}] {\it a timelike rotational surface defined by
 $$ \phi(u,v)= \(u\cosh v,\sqrt{1-u^{2}}\cosh f, u\sinh v, \sqrt{1-u^{2}}\sinh f\),$$
 where $u\in (0,1)$ and }
 $$f=\pm \int \frac{3du}{u^{1/3}(1-u^{2})\sqrt{9 u^{-2/3}-u^{2}+1}}.$$

\item[{\rm (9)}] {\it a timelike rotational surface defined by
 $$ \phi(u,v)= \( u\cosh v,\sqrt{u^{2}-1}\sinh f, u \sinh v,\sqrt{u^{2}-1}\cosh f\),$$
 where $u\in (1,\infty)$ and }
 $$f=\pm \int \frac{3du}{u^{1/3}(u^{2}-1)\sqrt{9 u^{-2/3}-u^{2}+1}}.$$

\item[{\rm (10)}] {\it a timelike rotational surface defined by
 $$ \phi(u,v)= \(\sqrt{1+u^{2}}\cos f, \sqrt{1+u^{2}}\sin f, u \cos v, u\sin v\),$$
 where $u\in (0,\infty)$ and }
 $$f=\pm \int \frac{3du}{u^{1/3}(1+u^{2})\sqrt{9 u^{-2/3}-u^{2}-1}}.$$

\item[{\rm (11)}]   {\it a timelike rotational surface defined by
 $$ \phi(u,v)= \(\frac{1}{2}\! \left\{u(v^{2}- f^{2})+\frac{1}{u}+u\right\}\! ,uf,\frac{1}{2}\!\left\{u(v^{2}- f^{2})+\frac{1}{u}-u\right\}\! ,uv \),$$
 where $u\in (0,3^{3/4})$ and }
 $$f= \int \frac{3du}{u^{2}\sqrt{9-u^{8/3}}}.$$
\end{enumerate}\end{theorem} 

\begin{remark} The surfaces given by  (3), (6) and (10) in Theorem \ref{S14.3} are spherical; the surfaces given by (1), (2), (5), (8) and (9) are hyperbolical, and  the surfaces given by (4), (7) and (11) are parabolic.
\end{remark}

\section{Biconservative hypersurfaces of Minkowski  4-space $\mathbb E^{4}_{1}$}\label{S15}

For biconservative hypersurfaces in Minkowski  4-space $\mathbb E^{4}_{1}$ with diagonalizable shape operator and two distinct principal curvatures, Fu and Turgay proved the following classification theorem in \cite{FT16}.
 
\begin{theorem}\label{T:15.1} {\it Let $M$ be a biconservative hypersurface of the Minkowski  4-space $\mathbb E^4_1$ with diagonalizable shape operator and two distinct principal curvatures. Then it is congruent to one of following hypersurfaces:}
\begin{enumerate}

\item[{\rm (1)}] $\mu(s,u,v)=\(\beta(s), s \cos u\sin v, s \sin u \sin v, s \cos v\);$

\item[{\rm (2)}] $\phi_2(s,u,v)=\(s\,\sin u\, {\sinh} v, s\, \sin u\,  {\cosh}v, s\cos u, \gamma(s)\);$

\item[{\rm (3)}]  $\phi_3(s,u,v)=\(s\, {\cosh}u ,s \sinh u\sin v, s  \sinh u\cos v, \delta(s)\);$

\item[{\rm (3)}]  $\phi_4(s,u,v)=\(\dfrac{s}{2} \(u^2+ v^2\)+s+\mu(s), s v, su,\dfrac{s}{2}  \( u^2+ v^2\)+\mu(s)\)$,
\end {enumerate}
{\it for some functions $\beta,\gamma,\delta,\mu$.}
\end{theorem}

The following example of  biconservative hypersurfaces in the Minkowski  4-space $\mathbb E^4_1$ with 3  distinct principal curvatures were constructed by Fu and  Turgay in \cite{FT16}.

\begin{example}\label{E:15.2} {\it Let $M$ be a hypersurface in $\mathbb E^4_1$ defined by
\begin{equation}\begin{aligned}\notag
\phi(s,u,v)=&\(\frac{1}{2} s( u^{2} +v^2)+a v^{2}+s+f, su,( s+2a)  v,\right.
\\&\hskip.3in \left. \frac{1}{2}  s( u^2+ v^2)+a v^2+f\), \quad \mathbb R\ni a\neq 0.
\end{aligned}\end{equation}
Then $M$ is  biconservative if and only if  either $M$ is spacelike and
\begin{equation*}
f=c_1 \left(\ln (s+2 a)- \ln s-\frac{a}{s}-\frac{a}{s+2 a}\right)-\frac{s}{2}
\end{equation*}
or it is timelike and
\begin{equation}\notag
f=c\int\limits_{s_0}^s\left(x(x+2a)\right)^{2/3}dx-\frac{s}{2},
\end{equation}
where $c\ne 0, s_{0}$ are real numbers.}

\end{example}

Fu and  Turgay  also classified in \cite{FT16} biconservative hypersurfaces in Minkowski  4-space $\mathbb E^4_1$ with diagonalizable shape operator and three distinct principal curvatures as follows.

\begin{theorem}\label{T:15.3} {\it
If  $M$ is a hypersurface of $\mathbb E^4_1$ with diagonalizable shape operator and three distinct principal curvatures, then $M$ is biconservative if and only if it is congruent to one of seven hypersurfaces:}
\begin{enumerate}
\item[(1)] {\it a generalized cylinder $N^2\times \mathbb E^1_1$ where $N^2$ is a biconservative  surface in $\mathbb E^3$;}

\item[(2)] {\it a generalized cylinder $N^2\times \mathbb E^1$ where $N^2$ is a biconservative spacelike surface in $\mathbb E^3_1$;}

\item[(3)] {\it a generalized cylinder $N^2_1\times \mathbb E^1$, where $N^2_1$ is a biconservative timelike surface in $\mathbb E^3_1$;}

\item[(4)] {\it a spacelike surface defined by
\begin{equation*}
\phi(s,u,v)=\left(s\,{\cosh}\, u ,s\,{\sinh}\, u ,f(s) \cos v,f(s)\sin v\right)
\end{equation*}
 for a function $f$ satisfying}
\begin{equation*}\frac{f''}{f'^2-1}=\frac{ff'+s}{sf};\end{equation*}

\item[(5)] {\it a timelike surface defined by
\begin{equation*}
\phi(s,u,v)=\left(s\,{\cosh}\, u ,s\,{\sinh}\, u ,f(s) \cos v,f(s)\sin v\right)\end{equation*}
 for a function $f$ satisfying}
\begin{equation*}\frac{-3f''}{f'^2-1}=\frac{ff'+s}{sf};\end{equation*}

\item[(6)] {\it a spacelike surface defined by
\begin{equation*}
\phi(s,u,v)=\left(s\,\sinh \,u,s\,{\cosh}\,u, f(s) \cos v, f_s)\sin v\right)
\end{equation*}
 for a function $f_3$ satisfying}
$$\frac{f''}{f'^2+1}=\frac{ff'+s}{sf};$$

\item[(7)] {\it a surface given in Example \ref{E:15.2}.}
\end{enumerate}
\end{theorem}

In \cite{KT23},  Kayhan and . Turgay constructed  an example of biconservative hypersurfaces in $\mathbb E^{4}_{1}$ with non-diagonalizable shape operator as follows.

\begin{proposition}\label{P:15.4} {\it Let $\gamma:I\to  \mathbb E^{4}_{1}$ be a null curve  which admits a pseudo-orthonormal frame field $\{T,U; v_{1}, v_{2}\}$ of vector fields defined along $\gamma$ such that
$$\<U,U\>=\<T,T\>=0,\;\;\; \<T,U\>=-1,$$
$$ \gamma'
=T,\;\;\; v'_{1}=\mu_{3} T,\;\;\; v'_{2}=\mu_{2}U+\mu_{4}T,$$
where $I$ is an open interval and $\mu_{2},\mu_{3},\mu_{4}: I\to \mathbb R$ are functions.
Then, the hypersurface $M$ in $\mathbb E^{4}_{1}$ 
parametrized by
$$\phi(s,u,w)=\gamma(s)+uT(s) +w v_{1}(s) +f(w)v_{2}(s), \;\; w \in I,\; (s,u)\in D,$$
has non-diagonalizable shape operator and it is proper biconservative if and
only if $f$ is a  function satisfying
$ 3ff''-2f'^{2}=2,$
where $D$ is an open subset of $\mathbb R^{2}$.}
\end{proposition}

The next classification theorem was due to  Kayhan and  Turgay  \cite{KT23}.

\begin{theorem}\label{T:15.5} {\it If $M$ is a hypersurface in $\mathbb E^{4}_{1}$, then $M$ has non-diagonalizable shape operator and it is proper biconservative if and only
if it is locally congruent to the hypersurface described by Proposition \ref{P:16.4}.}
\end{theorem}

\section{Quasi-minimal biconservative surfaces in $\mathbb E^4_2$}\label{S16}

The notion of {\it trapped surfaces} was introduced by  R. Penrose in \cite{Pe65}, which 
plays a important role in cosmology and general relativity. A black hole is a trapped region in a space-time enclosed by a marginally trapped surface.
 In term of mean curvature vector,  a spacelike surface in a given space-time is called {\it marginally trapped} if its mean curvature vector field is lightlike. More generally, a nondegenerate submanifold  of a pseudo-Riemannian manifold is called {\it quasi-minimal} if its mean curvature vector is lightlike at every point (see, e.g., \cite{Chen09,book11,CG09,R82}).

All flat biharmonic quasi-minimal surfaces in the  pseudo-Euclidean 4-space $\mathbb E^4_2$ endowed with a neutral metric were classified by the author in \cite{Chen08,CI98}. It follows from the proof of Theorem 4.1 of \cite{Chen08} that we have the following result (see \cite[Proposition 4.1]{YKTC21}).

\begin{theorem}\label{T:16.1} {\it A flat surface in $\mathbb E^4_2$ is quasi-minimal and biconservative   if and only if, locally, it is congruent to one of the following surfaces:}
\begin{enumerate}
\item[{\rm (1)}] {\it A surface defined by  $$\phi(s,t)=\(f(s,t),\frac{s-t}{\sqrt2},\frac{s+t}{\sqrt2},f(s,t)\),\ (s,t)\in U,$$  where $f:U\to\mathbb R$ is a function and $U$ is open in $\mathbb R^2$;}

\item[{\rm (2)}] {\it A surface given by  $$\phi(s,t)=\eta(s)t +\zeta(s),$$  where $\eta(s)$ is a lightlike curve in the light-cone $\mathcal LC$ and $\zeta$ is a light-like curve satisfying $\langle \eta',\zeta'\rangle=0$ and $\langle \eta,\zeta'\rangle=-1$.}
\end{enumerate}
\end{theorem}

In \cite{YKTC21},  Ye\u{g}in \c{S}en,  Kelleci,  Turgay and  Canfes provided the following existence result for quasi-minimal proper biconservative surfaces in $\mathbb E^4_2(c)$.

\begin{proposition}\label{P:16.2} {\it
Let $U=I_{1} \times I_{2}$ for open intervals $I_{1},I_{2}$, and $\beta\in C^\infty(U)$. Assume $\psi:U\to\mathbb R$ of $(U,g_\beta)$ satisfies $\psi(s,t)=\psi(t)$ for $(s,t)\in U$. Consider a light-like curve $\gamma:I_{2}\hookrightarrow \mathbb E^4_2$ lying on $\mathcal LC$ such that
$V_t={\rm Span}\{\gamma(t),\gamma'(t)\} $ is two dimensional for all $t\in I_{2}$. Suppose that  $\eta:I_{2}\to\mathbb R^4$ satisfies 
\begin{equation*}\langle \eta' ,\eta' \rangle=
\langle \gamma , \eta' \rangle=0,\end{equation*} 
\begin{equation*}\langle \eta' ,\gamma' \rangle=-\frac{1}{a }\;\; {\rm and}\;\;
\langle \eta' ,\gamma'' \rangle=\frac{2 a' -a  \psi}{a^2},
\end{equation*}
for a function $a\in C^\infty(I_{2})$. Then the map $\phi:(U,g_\beta)\to \mathbb E^4_2$ given by
\begin{equation}\notag
\phi(s,t)=\eta(t)+\big(s a'(t)-a(t) (\beta(s,t)+s \psi(t))\big)\gamma(t)+sa(t)\gamma'(t) 
\end{equation}
is a quasi-minimal, proper biconservative isometric immersion.}
\end{proposition}

The next result also from \cite{YKTC21} classified quasi-minimal, proper biconservative surfaces in $\mathbb E^4_2$.

\begin{theorem} \label{T:16.3} {\it  A  surface $M$ in $\mathbb E^4_2$ is quasi-minimal proper biconservative if and only if it is congruent to one of the following surfaces:}
\begin{enumerate}
\item[(1)] {\it A surface defined by  $$\phi(s,t)=\(f(s,t),\frac{s-t}{\sqrt2},\frac{s+t}{\sqrt2},f(s,t)\), \;\; (s,t)\in U\subset \mathbb R^{2},$$  where the function $f:U\to\mathbb R$ is defined on an open subset $U\subset \mathbb R^2$;}

\item[(2)] {\it A surface defined by  $$\phi(s,t)=\eta(s)t +\zeta(s),$$  where $\eta(s)$ is a lightlike curve in  light-cone $\mathcal LC$ and $\zeta$ is a lightlike curve satisfying $\langle \eta',\zeta'\rangle=0$ and $\langle \eta,w\zeta\rangle=-1$;}

\item[(3)] {\it A surface given in Proposition \ref{P:16.2}. }
\end{enumerate}
\end{theorem}

\section{Biconservative surfaces in $\mathbb E^{5}_{s}$}\label{S17}

Kayhan and Turgay \cite{KT22} studied hypersurfaces  in $\mathbb E^{5}_{1}$ with non-diagonalizable shape operator, and proved the following three non-existence results for such hypersurfaces.

\begin{theorem}\label{T:17.1} {\it There does not exist proper biconservative hypersurface in $\mathbb E^{5}_{1}$
 with the shape operator given by}
 $$A=\begin{pmatrix}
 -2\alpha & 1& 0 & 0 \\ 0& -2\alpha & 0 & 0
 \\ 0&0&-2\alpha& 0\\ 0&0& 0& 10\alpha
  \end{pmatrix},
 $$ {\it where $\alpha$ is the mean curvature  of the hypersurface.}
\end{theorem}

\begin{theorem}\label{T:17.2} {\it There does not exist  proper biconservative hypersurface in $\mathbb E^{5}_{1}$
 with the shape operator given by}
$$A=\begin{pmatrix} -2\alpha & 0& 0 & 0 \\ 0& -2\alpha & 1 & 0
 \\-1 &0&-2\alpha& 0\\ 0&0& 0& 10\alpha
  \end{pmatrix},
 $$ {\it where $\alpha$ is the mean curvature  of the hypersurface.}
\end{theorem}

\begin{theorem} \label{T:17.3} {\it There does not exist  proper biconservative hypersurface in $\mathbb E^{5}_{1}$
 with the shape operator given by}
$$A=\begin{pmatrix} 2\alpha & 0& 0 & 0 \\ 0&2\alpha & 1 & 0
 \\ -1 &0&2\alpha& 0\\ 0&0& 0& -2\alpha
  \end{pmatrix},
 $$ {\it where $\alpha$ is the mean curvature  of the hypersurface.}\end{theorem}

In  \cite{UT16}  Upadhyay and Turgay investigated  biconservative hypersurfaces of index 2 in $\mathbb E^{5}_{2}$.
They provided complete classification of biconservative hypersurfaces with diagonalizable shape operator at exactly three distinct principal curvatures. 

\begin{theorem}\label{T:17.4} {\it If $M^{4}_{2}$ is a biconservative hypersurface of index 2 in $\mathbb E^5_2$ with the shape operator given by 
$$S=\mathrm{diag}(\kappa_1,0,0,\kappa_4),\;\; \kappa_4\neq 0,$$
then it is congruent to one of the following 8 type of generalized cylinders over surfaces for some  functions $\mu=\mu(s)$ and $\rho=\rho(s)$:}
\begin{enumerate}
\item[{\rm (1)}]  $\phi(s, t, u, v) = (t, u,\mu \cos v,\mu\sin v,\rho),\;\; \mu'^{2}+\rho'^2=1$;
\item[{\rm (2)}] $\phi(s, t, u, v) = (\mu \sinh v, t, u,\mu\cosh v,\rho),\;\; \mu'^2+\rho'^2=1$;
\item[{\rm (3)}] $\phi (s, t, u, v) = (\rho, t, u, \mu \cos v, \mu \sin v), \;\; \mu'^2 - \rho'^2=-1$;

\item[{\rm (4)}] $\phi(s, t, u, v) = (\mu \cosh v, t, u,\mu \sinh v,\rho),\;\; \mu'^2-\rho'^2=1$;
\item[{\rm (5)}] $\phi(s, t, u, v) = (\mu \cos v,\mu \sin v,t,u,\rho ),\;\; \mu'^2-\rho'^2=1$;
\item[{\rm (6)}]\ $\phi(s, t, u, v) = (\mu \sinh v,\rho ,t,u,\mu \cosh v),\;\; \mu'^2-\rho'^2=-1$;
\item[{\rm (7)}] $\displaystyle \phi(s, t, u, v) = \left( \frac{v^2s}{2}+\rho+s,t,u,vs,\frac {v^2s}{2}+\rho\right),$ {\it such that $1-2\rho'<0$.}
\item[{\rm (8)}] $\phi(s, t, u, v) = \left(\dfrac{sv^2}{2}+\rho,sv,t,u,\dfrac{sv^2}{2}+\rho +s\right)$,  {\it such that $1+2\rho'<0$.}
\end{enumerate}
\end{theorem}

\begin{theorem}\label{T:17.5} {\it
Let $M$ be a hypersurface of index 2 in  $\mathbb E^5_2$. Assume that its shape operator has the form 
$$S=\mathrm{diag}(\kappa_1,\kappa_2,\kappa_2,0),\;\; \kappa_2\neq0.$$
Then, it is congruent to one of the following eight type of cylinders for some functions $\mu=\mu(s)$ and $\rho=\rho(s)$:}
\begin{enumerate}
\item[{\rm (1)}] $\phi(s,t,u,v)=(v, \mu \cosh t, \mu\sinh t \cos u, \mu\sinh t \sin u, \rho),\;\; \mu'^2-\rho'^2=1$;

\item[{\rm (2)}]$\phi(s,t,u,v)=(v,\rho, \mu\cos t, \mu\sin t \cos u,\mu \sin t \sin u),\;\; \mu'^2-\rho'^2=-1$;

\item[{\rm (3)}] $\phi(s,t,u,v)=(\mu\cosh t \sin u, \mu\cosh t \cos u, \mu\sinh t,\rho,v),\;\; \mu'^2-\rho'^2=1$;

\item[{\rm (4)}] $\phi(s,t,u,v)=(\rho, \mu\sinh t, \mu\cosh t\cos u, \mu\cosh t \sin u, v),\;\; \mu'^2-\rho'^2=-1$;

\item[{\rm (5)}] $\phi(s,t,u,v)=(v, \mu\sinh t, \mu\cosh t\cos u, \mu\cosh t \sin u, \rho),\;\; \mu'^2+\rho'^2=1$;

\item[{\rm (6)}] $\phi(s,t,u,v)=(\mu\sinh t \cos u, \mu\sinh t \sin u, \mu\cosh u, \rho, v),\;\; \mu'^2+\rho'^2=1$;

\item[{\rm (7)}]$\phi(s, t, u, v) = \left(\dfrac{s (t^2+u^2)}{2}+\rho ,v,st,su,\dfrac{s (t^2+u^2)}{2}+\rho-s\right),$ {\it where $\rho$ satisfies $1-2\rho'<0$;}

\item[{\rm (8)}]$\phi(s, t, u, v) = \left(\dfrac{s (t^2-u^2)}{2}+\rho ,st,su,v,\dfrac{s (t^2-u^2)}{2}+\rho+s\right),$ {\it where $\rho$ satisfies $1+2\rho'<0$;}
 \end{enumerate}
\end{theorem}

\begin{theorem}\label{T:17.6} {\it
Let $M$ be a hypersurface of index 2 in  $\mathbb E^5_2$. Assume that its shape operator has the form 
$$S=\mathrm{diag}(\kappa_1,\kappa_2,\kappa_2\kappa_4),\;\; \kappa_4\neq \kappa_2$$
for some non-vanishing smooth functions $k_1,k_2,k_4$.
Then, it is congruent to one of the following eight type of hypersurfaces for some functions $\mu=\mu(s)$ and $\rho=\rho(s)$:}
\begin{enumerate}
\item[{\rm (1)}] $\phi=(\rho\sinh v, \mu \cosh t, \mu \sinh t \cos u, \mu \sinh t \sin u, \rho\cosh v),\; \mu'^2-\rho'^2=1;$

\item[{\rm (2)}]  $\phi=\left(\rho\cos v, \phi_{2}\sin v, \mu\cos t, \mu \sin t \cos u, \mu \sin t \sin u\right),\; \mu'^2-\rho'^2=-1;$

\item[{\rm (3)}] $x=\left(\mu \cosh t \sin u, \mu \cosh t \cos u, \mu \sinh t, \rho\cos v, \rho\sin v\right),\; \mu'^2-\rho'^2=1;$

\item[{\rm (4)}]$\phi=\left(\rho\sinh v, \mu\sinh t, \mu\cosh t\cos u, \mu\cosh t \sin u, \rho\cosh v\right),\; \mu'^2+\rho'^2=1;$

\item[{\rm (5)}] $\phi=\! \left(\rho\cosh v, \mu\sinh t, \mu\cosh t\cos u, \mu\cosh t \sin u, \rho\sinh v\right)\!,\, \mu'^2\!-\!\rho'^2=-1;$

\item[{\rm (6)}] $\phi=\left(\mu\sinh t \cos u, \mu\sinh t \sin u, \mu\cosh u, \rho\cos v,\rho\sin v\right),\; \mu'^2+\rho'^2=1;$

\item[{\rm (7)}] {\it A hypersurface given by
\begin{equation}\begin{aligned}\notag
\phi(s,t,u,v)=&\left(\frac{ s }{2}\left( t^2+ u^2- v^2\right)-a v^2 + f, v (2a+ s), s  t, s  u,\right.\\
&\ \ \left.\frac{ s}{2} \left( t^2+ u^2- v^2 \right)-a v^2 +f -s
\right)
\end{aligned}\end{equation}
for a  constant $a\ne 0$ and a function $f=f(s)$ satisfying  $1-2f'<0$};

\item[{\rm (8)}] {\it A hypersurface given by
\begin{equation}\begin{aligned}\notag
\phi(s,t,u,v)=& \left(\frac{s}{2}\left( t^2- u^2- v^2\right) +a  v^2+f, s  t, s  u, v ( s-2a),\right.\\
&\ \;\;\;  \left.\frac{s}{2}\left( t^2- u^2- v^2\right) +a  v^2+f+s\right) 
\end{aligned}\end{equation}
for a constant $a\ne 0$ and a  function $f=f(s)$ satisfying $1+2 f'<0$.}
\end{enumerate}
\end{theorem}

\begin{remark} Upadhyay and Turgay \cite{UT16} also provided an explicit example of biconservative hypersurfaces in $\mathbb E^{n+1}_{2}$ with $n$ distinct principal curvatures.
\end{remark}

\section{Biconservative submanifolds in $\mathbb E^{m}_{1}$}\label{S18}

 In 2017, Deepika \cite{De17} studied biconservative Lorentzian hypersurface of $\mathbb E^{n+1}_{1}$ having non-diagonalizable shape operator with complex eigenvalues and proved the following two results.
 
 \begin{theorem}\label{T:18.1} {\it Let $M_{1}^{n}$ be a biconservative timelike
hypersurface in  $\mathbb E^{n+1}_{1}$ having non-diagonalizable shape operator with complex eigenvalues and with at most 5 distinct principal curvatures. Then $M_{1}^{n}$ has constant mean curvature.}
\end{theorem}

\begin{theorem}\label{T:18.2} {\it If $M_{1}^{n}$ is a biconservative  timelike  hypersurface  in $\mathbb E^{n+1}_{1}$ with constant length of second fundamental form and whose shape operator has complex eigenvalues with 6 distinct principal curvatures, then $M_{1}^{n}$  has constant mean curvature.}
\end{theorem}

In 2019,   Gupta and  Sharfuddin \cite{GS19} improved Deepika's results to the following.

\begin{theorem} \label{T:18.3} {\it If $M_{1}^{n}$ is a biconservative  timelike  hypersurface  in
$\mathbb{E}_{1}^{n+1}$ with complex eigenvalues, then it has constant mean curvature.} \end{theorem}

In  \cite{GS19}, Gupta and  Sharfuddin also provided some examples of such hypersurfaces.

\section{Biconservative PNMCV submanifolds in $\mathbb E^{n+2}_{1}$}\label{S19}

In \cite{Ka23},  Kayhan studied PNMCV submanifold $M^{n}$ with non-diagonalizable shape operator in $\mathbb E^{n+2}_{1}$ and proved the following.

\begin{theorem}\label{T:19.1}  {\it There is no PNMCV  biconservative submanifold $M^{n}$ in  $\mathbb E^{n+2}_{1}$ with non-diagonalizable shape operator such that the gradient of the mean curvature function is lightlike.}
\end{theorem}

Very recently,   Ye\u{g}in \c{S}en  \cite{S24} investigated  biconservative PNMCV spacelike submanifolds $M^{3}$ in $\mathbb E^{5}_{1}$ and established explicit classifications of these submanifolds with exactly two distinct principal curvatures of shape operator  in the direction of the mean curvature vector $H$.

\section{Biconservative submanifolds  in complex space forms}\label{S20}

A submanifold $M$ of a K\"ahler manifold $(\tilde M,J,g)$  equipped with the almost complex structure $J$ and K\"ahlerian metric $g$ is called \emph{totally real} if we have (see \cite{CO74}) $$J(T_{p}M)\subset T^{\perp}_{p}M,\;\; \forall p\in M.$$
In particular, if  $J(T_{p}M)=T^{\perp}_{p}M$, $\forall p\in M$, then $M$ is called {\it Lagrangian}.

A K\"ahler manifold $(\tilde M,J,g)$ is called a {\it complex space form} if it has constant holomorphic sectional curvature. The curvature tensor field 
$\tilde{R}$ of a complex space form $\tilde M(c)$ of constant holomorphic sectional curvature $c$ satisfies
\begin{equation}\begin{aligned} \tilde{R}(\bar{X},\bar{Y})\bar{Z} =&\,\frac{c}{4}\Big\{\langle \bar{Y},\bar{Z}\rangle \bar{X}-\langle \bar{X},\bar{Z}\rangle\bar{Y} +\langle J\bar{Y},\bar{Z}\rangle J \bar{X}\nonumber
\\
&\hskip.2in  \quad -\langle J \bar{X},\bar{Z}\rangle J\bar{Y} +2\langle J\bar{Y},\bar{X}\rangle J\bar{Z}\Big\},
\end{aligned}\end{equation}
for  vector fields $\bar{X}, \bar{Y},  \bar{Z} \in {\mathfrak X}(T\tilde M(c))$.

For a submanifold $M^{n}$ in a K\"ahler manifold, we put $$JH=(JH)^{T}+(JH)^{\perp},$$ where $(JH)^{T}$ and $(JH)^{\perp}$ denote the tangential and normal components of $JH$, where $H$ is the mean curvature vector.

Next,  we present the results on  biconservative submanifolds in complex space forms obtained by H. Bibi, the author, D. Fetcu and C. Oniciuc in \cite{BCFO23}.

\subsection{PMCV submanifolds  in complex $m$-space forms}\label{S20.1}

The following results were proved in \cite{BCFO23}.

\begin{theorem}\label{T:20.1} {\it If  $M^{n}$ is a PMCV submanifold of a complex space form $\tilde M^{m}(c)$, then}
\begin{itemize}
  \item [{\rm (1)}] {\it if $c=0$,  $M^{n}$ is a biconservative submanifold;}
    \item [{\rm (2)}] {\it if $c\neq 0$,  $M^{n}$ is biconservative if and only if $J(JH)^{T} \in  {\mathfrak X}(T^{\perp}M^{n})$.}
\end{itemize}\end{theorem}

Theorem \ref{T:20.1} implies the following.

\begin{corollary} \label{C:20.2} {\it  If $M^n$ be a PMCV totally real submanifold of a complex space form, then $M^n$ is biconservative.}
\end{corollary}

\begin{corollary} \label{C:20.3} {\it  Every PMCV real hypersurface $M^{2n-1}$ of a complex space form  $\tilde M^{n}(c)$ is biconservative.}
\end{corollary}

\begin{theorem} \label{T:20.4} {\it  Let $M^{n}$ be a PMCV submanifold of  a complex space form $\tilde M^{m}(c)$ with $c \neq 0$. If $JH \in {\mathfrak X}(T^{\perp}M^{n})\, ($or if $JH \in  {\mathfrak X}(TM^{n}))$, then $M^{n}$ is biconservative.}
\end{theorem}

\begin{theorem}\label{T:20.5} {\it Let $M^{2}$  be a PMCV surface in  a complex space form  $\tilde M^{m}(c)$. Then}

\begin{itemize}
  \item [{\rm (1)}] {\it If $c=0$,  $M^{2}$ is biconservative.}
\item [{\rm (2)}] {\it If $c\neq 0$, then $M^{2}$ is biconservative if and only if $M^{2}$ is totally real.}
\end{itemize}\end{theorem}

Given a submanifold $M$ of a K\"ahler manifold $\tilde M$,  for a  vector $0\ne X$ of $M$ at $p$, the angle $\theta (X)$  between
$JX$ and $T_pM$ is called the {\it Wirtinger angle} of $X$.
The submanifold $M$ is called a {\it slant submanifold\/} if the Wirtinger angle $\theta (X)$ is independent of the choice of $X \in T_{p}M$ and of $p \in M$ (see \cite{chen90,book90}.)

For slant surfaces, Theorem \ref{T:20.5} implies the following.

\begin{corollary} \label{C:20.6} {\it  Every PMCV proper slant surface in a non-flat complex space form  is never biconservative.}
\end{corollary}

\subsection{$CMC$ biconservative surfaces in complex 2-space forms}\label{S20.2}

For $CMC$ biconservative surfaces in a complex 2-space form $\tilde M^{2}(c)$, we have the following results also from \cite{BCFO23}.

\begin{proposition} \label{P:20.7} {\it 
Let $M^{2}$ be a pseudo-umbilical CMCV biconservative  surface in $\tilde M^{2}(c)$ with $c\neq 0$. Then $M^{2}$ is PMCV and $J(JH)^{T}$ is normal.}
\end{proposition}

\begin{proposition}  \label{P:20.8} {\it  Let $M^{2}$ be a CMC biconservative surface with no pseudo-umbilical points in $\tilde M^{2}(c)$,  $c\neq 0$. If $J(JH)^{T}$ is normal, then $M^{2}$ is PMCV.}
\end{proposition}

\begin{theorem}\label{T:20.9}  {\it  Let $M^{2}$ be a CMC biconservative surface in $\tilde M^{2}(c)$. If $J(JH)^{T}$ is normal, then $M^{2}$ is both PMCV and totally real.}
\end{theorem}

For CMC but non-PMCV biconservative surfaces in $\mathbb C^{2}$, we have 

\begin{proposition} \label{P:20.10} {\it  If $M^{2}$  is  CMC,  non-PMCV biconservative surface  in  $\mathbb{C}^{2}$, then $(J(JH)^{T})^{T}\neq 0$.}
\end{proposition}

\subsection{Reduction of codimension for biconservative surfaces in complex space forms}\label{S20.3}

The following reduction theorem for biconservative surfaces  in complex space forms was also proved in \cite{BCFO23}.

\begin{theorem}\label{T:20.11} {\it If $M^{2}$ is a non-pseudoumbilical, PMCV, totally real surface in $ \tilde M^{n}(c)$ with $c\neq 0$ and $n\geq 4$, then there exists a totally geodesic complex submanifold $ \tilde M^{4}(c)\subset  \tilde M^{n}(c)$ such that $M^{2}\subset \tilde M^{4}(c)$.  }
\end{theorem}

We can improve Theorem \ref{T:20.11} and reduce the codimension even more under a slightly stronger assumption as follows.

\begin{theorem}\label{T:20.12} {\it
Let $M^{2}$ be a non-pseudoumbilical, PMCV, totally real surface in $\tilde M^{n}(c)$, $c\neq 0$. If $H \in {\mathfrak X}(J(TM^{2}))$, then there is a totally geodesic complex submanifold $ \tilde M^{2}(c)\subset  \tilde M^{n}(c)$ such that $M^{2}\subset \tilde M^{2}(c)$. } \end{theorem}

\subsection{Examples of non-PMC, CMC biconservative submanifolds in $CP^{m}$}\label{S20.4}

The Segre embedding of two complex projective spaces $CP^{h}(4)$ and $CP^{p}(4)$:
\begin{equation*} \label{1.1} S_{hp}:CP^{h}(4)\times CP^{p}(4)\to CP^{h+p+hp}(4),\end{equation*}
is defined by 
\begin{equation*} \label{1.2} S_{hp}(z_0,\ldots,z_{h};w_0,\ldots,w_p)=\big(z_jw_t \big)_{0\leq j\leq h,0\leq t\leq p},\end{equation*} 
where  $(z_0,\ldots,z_{h})$ and
$(w_0,\ldots,w_p)$ are the homogeneous coordinates of
$CP^{h}(4)$ and $CP^p(4)$, respectively. 
This embedding was introduced by C.  Segre in 1891 (see \cite{Se1891}). It is well-known  that the Segre embedding $S_{hp}$ is also a K\"ahlerian
embedding. 
When $h=p=1$,
the Segre embedding is nothing but the complex quadric
surface,
$Q_2=CP^1\times CP^1$ in $CP^3$, defined by
\begin{equation*} 
Q_2=\Bigg\{(z_0,z_1,z_2,z_3)\in CP^3: \sum_{j=0}^3
z_j^2=0\Bigg\}.\end{equation*}

Bibi et. al. constructed in \cite{BCFO23} examples with the Segre embedding  of non-PMC, CMC, biconservative submanifolds in complex projective spaces as follows.

\begin{theorem} \label{T:20.13} {\it
Let $\gamma$ be a curve of nonzero curvature $\kappa$ in $CP^1(4)$. Then we have:}
\begin{enumerate}
\item[{\rm (a)}] {\it via the Segre embedding $S_{1p}$ of $CP^{1}(4)\times CP^{p}(4)$ into 
$\mathbb{C}P^{1+2q}(4)$, the image of $M^{1+2q}=\gamma \times CP^p (4)$ under $S_{1p}$ is a  biconservative submanifold of 
$CP^{1+2p}(4)$ 
if and only if $\kappa$ is constant. Moreover,  $M^{1+2p}$ is CMC, non-PMC, and  it is not totally real;}
\item[{\rm (b)}] {\it $M^{1+2p}$ is a  proper biharmonic submanifold of 
$CP^{1+2p}(4)$ if and only if 
$\kappa^{2}=4$, i.e., $\gamma$ is proper-biharmonic in $\mathbb{C}P^{1}(4)$}
\end{enumerate}  
\end{theorem}

These examples show that the existence of biconservative submanifolds of dimension greater than two is much less rigid.

\section{$L_{k}$-biconservative hypersurfaces in $\mathbb E^{n+1}_{1}$}\label{S20}

\subsection{Shape operator of timelike hypersurface $M^{3}_{1}$}\label{S21.1}

The shape operator of the timelike hypersurface $M^{3}_{1}$ in $\mathbb E^{4}_{1}$ can be put into one of the four matrix forms,  denoted by Types I, II, III and IV. 

For Type I and Type IV,  the matrix representation of the induced metric on $M^{3}_{1}$, with respect to an orthonormal frame of the tangent bundle $TM^{3}_{1}$, is
$$\begin{pmatrix} -1&0&0\\ 0&1&0\\0&0&1\end{pmatrix}$$ 
and the corresponding shape operator $S$ can be put into matrix forms as
$$S_{1}=\begin{pmatrix} \lambda_{1}&0&0\\ 0&\lambda_{2}&0\\ 0&0&\lambda_{3}\end{pmatrix},
\;\;\; S_{4}=\begin{pmatrix} \kappa&\mu&0\\ -\mu &\kappa&0\\0&0&\eta\end{pmatrix},\;\; (\mu\ne 0),$$ respectively.
For Type II and Type III,  the matrix representation of the induced metric
on $M^{3}_{1}$, with respect to a pseudo-orthonormal frame of 
 of $TM^{3}_{1}$, is
$$\begin{pmatrix} 0&1& 0 \\ 1&0&0\\0&0&1\end{pmatrix}$$ 
and the corresponding  shape operator $S$ can be put into matrix forms as
$$S_{2}=\begin{pmatrix} \kappa&0&0\\ 1&\kappa&0\\ 0&0&\mu \end{pmatrix},
\;\;\; S_{3}=\begin{pmatrix} \kappa&0&0\\ 0&\kappa&1\\-1&0&\kappa\end{pmatrix},$$ respectively.

\subsection{$k$-th mean curvature of hypersurfaces $M^{n}_{1}$}\label{S21.2}
Following the four possible matrix representations of the shape
operator $S_{i}\, (i=1,2,3,4)$ of $M^{3}_{1}$ in the previous subsection, we define its principal curvatures, denoted by unified notations
$\kappa_{i}$ for $i=1,2,3$, according to \cite{Pa22}, as follow.

\vskip.06in
For Type I, we put $\kappa_{i}=\lambda_{i}$, $i=1,2,3$,  where $\lambda_{i}$'s were given in $S_{1}$.

For Type II, we put  $\kappa_{1}=\kappa_{2}=\kappa$ and $\kappa_{3}=\mu$, where $\mu$ was given in $S_{2}$.

For Type III, we put  $\kappa_{1}=\kappa_{2}=\kappa_{3}=\kappa$, where $\kappa$ was given in $S_{3}$.

For Type IV, we put $\kappa_{1}=\kappa_{1}+i\mu$, $\kappa_{2}=\kappa_{1}-i\mu$
 and $\kappa_{3}=\eta$, where $\kappa,\mu,\eta$ were given in $S_{4}$.
\vskip.06in

The {\it characteristic polynomial} of the shape operator $S$ on $M^{3}_{1}$
 is 
$$Q(t)=\prod_{i=1}^{3} (t-\kappa_i)=\sum_{k=0}^{3}(-1)^{k}s_{k}t^{3-k},$$
where 
$$s_{0=1},\;\; s_{1}=\sum_{j=1}^{3}\kappa_{j},\;\; s_{2}=\sum_{i\leq i_{1}\leq i_{2}\leq 3}\kappa_{i_1\kappa_{i_2}}\;\; {\rm and}\;\; s_{3}=\kappa_{1}\kappa_{2}\kappa_{3}.$$ 
The $k$-th mean curvature of $M^{3}_{1}$ is given by $H_{k}=s_{k}/\binom{3}{k}$,  $k\in \{1,2,3\}$. 
$M^{3}_{1}$ is called {\it $(k-1)$-minimal} if $H_{k}=0$ holds identically.

Similar definitions holds for hypersurfaces $M^{n}_{1}$ with $n>3$.

\subsection{Definition of $L_{k}$-biconservative  hypersurfaces}\label{S21.3}

The operator $L_{k}$ is an extension of the ordinary Laplace operator $\Delta \, (=L_{0})$ which stands for the linearized operator of the first variation of $(k+1)$-th mean curvature function  (see, e.g.,  \cite{book15,Re73}). 

The operator $L_{k}$ is defined on $M^{n}_{1}$  by 
$$L_{k}(f) = {\rm Trace}(P_{k}\circ \nabla^{2}f),\;\;\; \in C^{\infty}(M^{n}_{1}),$$
where $P_{k}$  is the $k$-th Newton transformation associated with the second fundamental from of the hypersurface and $\nabla^{2}f$ is the Hessian of $f$ (see \cite{book15,Re73} for details). 

The next definition extended the notions of biharmonic and biconservative hypersurfaces to $L_{k}$-biharmonic and $L_{k}$-biconservative hypersurfaces.

\begin{definition} A hypersurface $\phi: M^{n}_{1}\to \mathbb E^{n+1}_{1}$
is called {\it $L_{k}$-biharmonic} if it satisfies  $L^{2}_{k}\phi=0$; and $\phi$ is called 
{\it $L_{k}$-biconservative} if the tangent component of $L^{2}_{k}\phi$ vanishes identically. 
\end{definition}

\subsection{$L_{k}$-biconservative timelike hypersurfaces in $\mathbb E^{4}_{1}$}\label{S21.4}

In \cite{Pa22}, F. Pashaie proved the following results.

\begin{theorem} \label{T:21.2}
 {\it Let $k\in \{0,1,2\}$ and let $\phi: M^{3}_{1}\to \mathbb E^{4}_{1}$ be a $L_{k}$-biconservative hypersurface in 
 $\mathbb E^{4}_{1}$, having diagonalizable shape operator $($i.e.,
of type I\/$)$ with constant $k$-th mean curvature and exactly two distinct principal
curvatures. Then, it
has constant $(k + 1)$-th mean curvature.}\end{theorem}

\begin{theorem} \label{T:21.3} {\it Let $k\in \{0,1,2\}$ and let 
$\phi: M^{3}_{1}\to \mathbb E^{4}_{1}$ be an $L_{k}$-biconservative  hypersurface
with type II shape operator  which has exactly two distinct
principal curvatures and constant $k$-th mean curvature, then its $(k + 1)$-th
mean curvature is constant.}
\end{theorem}

\begin{theorem} \label{T:21.4} {\it Let $k\in \{0,1,2\}$ and let 
$\phi: M^{3}_{1}\to \mathbb E^{4}_{1}$ be an $L_{k}$-biconservative timelike hypersurface
with  type III  shape operator which has constant $k$-th mean
curvature, then its $(k + 1)$-th mean curvature is constant.}
\end{theorem}

\begin{theorem}\label{T:21.5} {\it Let $k\in \{0,1,2\}$ and let $\phi: M^{3}_{1}\to \mathbb E^{4}_{1}$ be an $L_{k}$-biconservative timelike hypersurface with type IV shape operator which has constant $k$-th mean curvature and a constant real principal curvature. Then, its second and third mean curvatures are constant.}
\end{theorem}

\subsection{$L_{k}$-biconservative timelike hypersurfaces in $\mathbb E^{5}_{1}$}\label{S21.5}

Anagous to hypersurface $M^{3}_{1}\subset \mathbb E^{4}_{1}$, shape operators of hypersurfaces $M^{4}_{1}$ in $\mathbb E^{5}_{1}$ can also be divided into types I, II, III and IV.

For $L_{k}$-biconservative  hypersurfaces in $\mathbb E^{5}_{1}$, Pashaie proved the following results in \cite{Pa20}.
 
\begin{theorem} \label{T:21.6} 
 {\it Let $\phi: M^{4}_{1}\to \mathbb E^{5}_{1}$ be a $L_{k}$-biconservative orientable hypersurface $($for a nature number  $k \leq 3)$ in $\mathbb E^{5}_{1}$ having diagonalizable shape operator.
If $M^{4}_{1}$  has a principal curvature of multiplicity four, then it has constant
$(k + 1)$-th mean curvature.}\end{theorem}

\begin{theorem} \label{T:21.7} 
 {\it Let $\phi: M^{4}_{1}\to \mathbb E^{5}_{1}$ be a $L_{k}$-biconservative
orientable  hypersurface $($for a nature number  $k \leq 3)$ in  $\mathbb E^{5}_{1}$
having constant ordinary mean
curvature and diagonalizable shape operator. If $M^{4}_{1}$ has two principal curvatures
 both of multiplicity $2$, then it has constant $(k+1)$-th
mean curvature.}\end{theorem}

\begin{theorem} \label{T:21.8} 
 {\it Let $\phi: M^{4}_{1}\to \mathbb E^{5}_{1}$ be a $L_{k}$-biconservative
 Lorentzian hypersurface $($for a nature number  $k \leq 3)$ in  $\mathbb E^{5}_{1}$
having diagonalizable shape operator
with constant ordinary mean curvature. If $M^{4}_{1}$ has exactly two principal
curvatures  of multiplicities $3$ and $1$ $($respectively$)$, then it has constant $(k+1)$-th mean curvature.}\end{theorem}

\begin{theorem} \label{T:21.9} 
 {\it Let $\phi: M^{4}_{1}\to \mathbb E^{5}_{1}$ be a $L_{k}$-biconservative
 hypersurface $($for a nature number  $k \leq 3)$ in   $\mathbb E^{5}_{1}$
having type II shape operator.
If $M^{4}_{1}$  has at most two distinct principal curvatures and constant ordinary
mean curvature, then it has constant $(k+1)$-th mean curvature.}
 \end{theorem}

\begin{theorem} \label{T:21.10} 
 {\it Let $\phi: M^{4}_{1}\to \mathbb E^{5}_{1}$ be a $L_{k}$-biconservative
hypersurface $($for a nature number   $k \leq 3)$ in  $\mathbb E^{5}_{1}$
having  type III shape operator. If
$M^{4}_{1}$ has constant ordinary mean curvature, then it has constant $(k+1)$-th
mean curvature.}
\end{theorem}

\begin{theorem} \label{T:21.11} 
 {\it Let $\phi: M^{4}_{1}\to \mathbb E^{5}_{1}$ be a $L_{k}$-biconservative
hypersurface $($for a nature number  $k\leq 3)$ in  $\mathbb E^{5}_{1}$
having  type IV shape operator. If $M^{4}_{1}$  has at most two distinct non-zero principal curvatures, then it has constant $(k+1)$-th
mean curvature and $M^{4}_{1}$  is isoparametric.}
\end{theorem}

\subsection{$L_{1}$-biconservative hypersurfaces in $\mathbb E^{4}_{1}$}\label{S21.6}

Very recently, F. Pashaie studied $L_{1}$-biconservative hypersurfaces in $\mathbb E^{4}_{1}$ and proved the following results in \cite{Pa23,Pa24}.

\begin{theorem} \label{T:21.12} {\rm \cite{Pa23}}   {\it Every $L_{1}$-biconservative hypersurface $M^{4}_{1}$ in $\mathbb E^{5}_{1}$  with constant mean curvature and at least three principal
curvatures has constant second mean curvature.}
\end{theorem} 

\begin{theorem} \label{T:21.13} {\rm \cite{Pa24}}   {\it  If $M$ is a $L_{1}$-biconservative spacelike hypersurface in $\mathbb E^{4}_{1}$ with constant  mean curvature and  at most two distinct principal curvatures, then the scalar curvature of $M$ is constant and $M$ is isoparametric.}
\end{theorem} 

\begin{theorem} \label{T:21.14} {\rm \cite{Pa24}}  {\it  If $M$ is a $L_{1}$-biconservative spacelike hypersurface in $\mathbb E^{4}_{1}$ with constant  mean curvature and three distinct principal curvatures, then the scalar curvature of $M$ is constant.}
\end{theorem}

\section{Generalized biconservative surfaces in $\mathbb E^{3}_{1}$}\label{S22}

Recall from Corollary \ref{C:4.10} that  a hypersurface of a real space form $R^{n+1}(c)$ is biconservative if and only if its mean curvature $\alpha$ satisfies $$A(\nabla \alpha)=-\frac{n}{2}\alpha\nabla \alpha.$$ 
In \cite{YZ21}, D. Yang and X. Y. Zhu defined a {\it generalized biconservative hypersurface} in $R^{n+1}(c)$ to be a hypersurface satisfying  $$A(\nabla)=k \alpha \nabla \alpha$$
for some constant $k$.

For generalized biconservative hypersurfaces, Yang and Zhao  \cite{YZ21} proved the following results.

\begin{theorem} \label{T:22.1}   {\it Let  $\phi: M\to \mathbb E^{3}_{1}$ be a non-CMC spacelike generalized biconservative surface in $\mathbb E^{3}_{1}$. Then $\phi$ is locally given by one of the following surfaces:}
\begin{enumerate}

\item[{\rm (1)}] {a cylinder defined by $$\phi(s,t)=\(t,\int^{s}\cosh f(s)ds, \int^{s}\sinh f(s)ds\),$$ where $f=f(s)$ is a function satisfying $f''\ne 0$};

\item[{\rm (2)}] {\it a cone given by
$$\phi(s,t)=\frac{a-bs}{\sqrt{4-b^{2}}}\(\frac{2}{b}, \sinh (\sqrt{4-b^{2}} t),\cosh (\sqrt{4-b^{2}} t)\),$$
where $a,b$ are constants and $b \in (-2, 2)$;}

\item[{\rm (3)}] {\it a cone defined by 
$$\phi(s,t)=\frac{a-bs}{\sqrt{b^{2}-4}}\( \cos (\sqrt{b^{2}-4} t),\sin (\sqrt{b^{2}-4} t),\frac{2}{b}\),$$
where $a,b$ are constants and $b\in \mathbb R\setminus [-2,2]$;}

\item[{\rm (4)}] {\it a surface of revolution with a lightlike axis defined by
$$\phi(s,t)=(c\pm 2s)\( \frac{1}{2} -t^{2}, t, t^{2}\),$$
where $c$ is any constant;}

\item[{\rm (5)}] {\it a surface of revolution with a timelike axis defined by
$$\phi(s,t)=(s\cos t, s\sin t, f(s)),$$
where 
$$f(s)=\int^{s} \frac{(2+k)ds} {\sqrt{s^{-2k/(k+2)}+(2+k)^{2}}}$$
and $k\in \mathbb R\setminus \{-2,-1,0\}$;}

\item[{\rm (6)}] {\it a surface of revolution with a spacelike axis defined by
$$\phi(s,t)=( f(s), s\sinh t, s\cosh t),$$
where 
$$f(s)=\int^{s} \frac{(2+k)ds} {\sqrt{(2+k)^{2}-s^{2k/(k+2)}}}$$
and $k \in \mathbb R\setminus \{-2,-1,0\}$;}

\item[{\rm (7)}] {\it a surface of revolution with a lightlike axis defined by
$$\phi(s,t)=\(\frac{1}{2}st^{2}-\frac{1}{p(k)}s^{q(k)}-\frac{s}{2}, st, \frac{1}{2}st^{2}-\frac{1}{p(k)}s^{q(k)}+\frac{s}{2}\),$$
where $p(k) =2(k+2)(3k+2)$, $q(k)=\frac{3k+2}{k+2}$, and $k\in \mathbb R\setminus \{0,1,2\}$.}

\end{enumerate}
\end{theorem}

\begin{theorem} \label{T:22.2}   {\it Let $\phi: M_{1}\to \mathbb E^{3}_{1}$ be a non-CMC timelike generalized biconservative surface in $\mathbb E^{3}_{1}$. Then $\phi$ is locally one of the following surfaces:}

\begin{itemize}
\item[{\rm (1)}] {\it a cylinder defined by
$$\phi(s, t) = \(\int^{s} \sin f(s)ds,\int^{s} \cos f(s)ds,t\),$$ where $f(s)$ satisfis $f''\ne 0$;}

\item[{\rm (2)}]{\it a cylinder defined  by $$\phi(s, t) = \(\int^{s} \sinh f(s)ds,t,\int^{s} \cosh f(s)ds\),$$ where $f(s)$ satisfies $f''\ne 0$;}

\item[{\rm (3)}] {\it a cone defined by
$$\phi(s,t)=\frac{a-bs}{\sqrt{b^{2}+4}}\(\frac{2}{b},\cosh (\sqrt{b^{2}+4}t),\sinh (\sqrt{b^{2}+4}t\)
$$
where $a,b$ are constants and $b\in \mathbb R$;}

\item[{\rm (4)}]  {\it a cone defined by
$$\phi(s,t)=\frac{a-bs}{\sqrt{b^{2}-4}}\(\frac{2}{b},\sinh (\sqrt{b^{2}-4}t),\cosh (\sqrt{b^{2}-4}t\)
$$
where $a,b\in \mathbb R$;}

\item[{\rm (5)}]  {\it a cone defined by
$$\phi(s,t)=\frac{a-bs}{\sqrt{4-b^{2}}}\(\cos (\sqrt{4-b^{2}}t),\sin (\sqrt{4-b^{2}}t,\frac{2}{b}\)
$$
where $a,b$ are constants and $b\in (-2,2)$;}

\item[{\rm (6)}] {\it a surface of revolution with a lightlike axis defined by
$$\phi(s,t)=(a \pm 2s)\(-t^{2}, t, t^{2}+\frac{1}{2}\),$$
where $a$ is a constant;}

\item[{\rm (7)}] {\it a surface of revolution with a spacelike axis defined by
$$\phi(s,t)=( f(s), s\cosh t, s\sinh t),$$
where }
$$f(s)=\int^{s} \frac{(2-k)ds} {\sqrt{s^{2k/(k+2)}-(2+k)^{2}}};$$

\item[{\rm (8)}] {\it a surface of revolution with a spacelike axis defined by
$$\phi(s,t)=(  f(s), s\sinh t, s\cosh t),$$
where }
$$f(s)=\int^{s} \frac{(2-k)ds} {\sqrt{s^{2k/(k+2)}+(2-k)^{2}}};$$

\item[{\rm (9)}] {\it a surface of revolution with a timelike axis defined by
$$\phi(s,t)=( s\cos t, s\sin t, f(s)),$$
where }
$$f(s)=\int^{s} \frac{(2-k)ds} {\sqrt{(2-k)^{2}-s^{2k/(k+2)}}};$$

\item[{\rm (10)}] {\it a surface of revolution with a lightlike axis defined by
$$\phi(s,t)=\(\frac{1}{2}st^{2}+\frac{s^{q(k)}}{p(k)}-\frac{s}{2}, st, \frac{1}{2}st^{2}+\frac{s^{q(k)}}{p(k)}+\frac{s}{2}\),$$
where $p(k) =2(k-2)(3k-2)$, $q(k) =\frac{3k-2}{k-2}$, and $k\in \mathbb R\setminus \{0,1,2\}$.}

\item[{\rm (11)}] {\it a null scroll with non-constant mean curvature.}
\end{itemize}\end{theorem}

\begin{remark} A {\it null scroll} in a Minkowski 3-space $\mathbb E^{3}_{1}$ is a nondegenerate ruled surfaces over a null curve in $\mathbb E^{3}_{1}$.
\end{remark}

\section{Generalized biconservative surfaces in $S^{3}_{1}(1)$ and $H^{3}_{1}(-1)$}\label{S23}

D. Yang and Z.  Zhao \cite{YZ23} investigated generalized biconservative spacelike surfaces in  Lorentzian 3-space forms.

\subsection{Generalized biconservative surfaces in de Sitter 3-space $S^{3}_{1}(1)$}\label{S23.1}

In \cite{YZ23}, Yang and Zhao proved the following classification theorem for 
 generalized biconservative spacelike surfaces in  $S^{3}_{1}(1)$. 

\begin{theorem} \label{T:23.1}   {\it If $\phi:M\to S_{1}^{3}(1)\subset \mathbb E^{4}_{1}$
is a non-CMC generalized biconservative spacelike surface in the
de Sitter space $S_{1}^{3}(1)$, then $\phi=\phi(s,t)$ is locally one of the following surfaces of revolution:}
\begin{itemize}
\item[{\rm (1)}] {\it $\phi=\frac{1}{\sqrt{c}}\(\sqrt{4\!- \! b^{2}}\cos s+\frac{ab}{\sqrt{4-b^{2}}}\sin s,f \sinh \sqrt{c}t, \frac{2\sqrt{c}}{\sqrt{4-b^{2}}}\sin s, f\cosh \sqrt{c}t\)$, where $c=4-a^{2}-b^{2} >0$, $f=a\cos s +b \sin s$, and $a,b$ are constants;}

\item[{\rm (2)}] {\it $\phi\,=\,\frac{1}{\sqrt{c}}\,\(\sqrt{ b^{2}- 4}\cos s-\! \frac{ab}{\sqrt{b^{2}-4}}\sin s,f \cos \sqrt{c}t,f\sin \sqrt{c}t, \frac{2\sqrt{c}}{\sqrt{b^{2}-4}}\sin s\)$, where $c=a^{2}+b^{2}-4$, $f=a\cos s +b \sin s$, and $a,b$ are constants  such that $b\in \mathbb R\setminus [-2,2]$;}

\item[{\rm (3)}] {\it $\phi\,=\, \frac{1}{\sqrt{c}}\(\frac{2\sqrt{c}}{\sqrt{4-b^{2}}}\sin s, f\cos \sqrt{c}t,  f\sin \sqrt{c}t, \sqrt{4-b^2}\cos s +\frac{ab}{\sqrt{4-b^{2}}}\sin s\)$,\ where $c=a^{2}+b^{2}-4>0$, $f=a\cos s +b \sin s$, and $a$ and $b$ are constants with $a\ne 0$ and $b\in \{-2,2\}$;}

\item[{\rm (4)}] {\it $\phi=\(\cos s-\dfrac{a}{4}\sin s, \dfrac{f}{a}  \cos a t,  \dfrac{f}{a} \sin a t, \cos s+\( \dfrac{2}{a}-  \dfrac{a}{4}\)\sin s\)$, where \\$f=a\cos s +b \sin s$, and $a,b$ are constants  such that $b\in (-2,2)$;}

\item[{\rm (5)}] {\it $\phi=\(\cos s-\dfrac{a }{2}f t^{2}, \sin s-\dfrac{b}{2} ft^{2}, ft, f t^{2}\)$, where $f=a\cos s +b \sin s$, and $a,b$ are constants  satisfying $a^{2}+b^{2}=4$;  }

\item[{\rm (6)}] {\it  $\phi=\(s\sinh t, \sqrt{1+s^2}  \cos f, \sqrt{1+s^2}  \sin f, s\cosh t\)$, where 
$$f=\pm \int^{s}\frac{(2+k)ds} {(1+s^2)s^{k/(2+k)}\sqrt{(2+k)^2 s^{-2k/(+k)}-s^2-1}},$$ with
$s\in (0,1)$ and $k\ne 0,-1,-2$;}

\item[{\rm (7)}] {\it  $\phi=\(s\cos t, s\sin t, \sqrt{1-s^2}  \cosh f, \sqrt{1-s^2}  \sinh f\)$, where 
$$f=\pm \int^{s}\frac{(2+k)ds} {(1-s^2)s^{k/(2+k)}\sqrt{(2+k)^2 s^{-2k/(2+k)}-s^2+1}},$$ with
$s\in (0,1)$ and $k\ne 0,-1,-2$;}

\item[{\rm (8)}] {\it  $\phi=\(s\cos t, s\sin t, \sqrt{s^2-1}  \sinh f, \sqrt{s^2-1}  \cosh f\)$, where 
$$f=\pm \int^{s}\frac{(2+k)ds} {(s^2-1)s^{k/(2+k)}\sqrt{(2+k)^2s^{-2k/(2+k)}-s^2+1}},$$ with
$s\in (1,\infty)$ and $k\ne 0,-1,-2$;}

\item[{\rm (9)}] {\it  $\phi=\(\frac{1}{2}\(st^{2}+sf^{2}-\frac{1}{s}-s\), st, sf,\frac{1}{2}\(st^{2}+sf^{2}-\frac{1}{s}+s\)\)$, where 
$$f= \int^{s}\frac{(2+k)ds} {s^2\sqrt{(2+k)^2-s^{(4+4k)/(2+k)}}},$$ 
$s\in (0,(2+k)^{(2+k)/(2+2k)})$, and $k\ne 0,-1,-2$.}
\end{itemize}\end{theorem}

\subsection{Generalized biconservative surfaces in anti-de Sitter 3-space $H^{3}_{-1}(1)$}\label{S23.2}

In \cite{YZ23}, Yang and Zhao also obtained the following classification theorem for 
 generalized biconservative spacelike surfaces in  $H^{3}_{1}(-1)$. 

\begin{theorem} \label{T:23.2}   {\it If $\phi:M\to H_{1}^{3}(-1)\subset \mathbb E^{4}_{2}$ is a non-CMC generalized biconservative spacelike surface in anti-de Sitter space $H_{1}^{3}(-1)$,
then $\phi$ is locally one of the following surfaces of revolution:}
\begin{itemize}
\item[{\rm (1)}] {\it $\phi=\frac{1}{\sqrt{c}}\! \(\! f\cos \sqrt{c} t, f\sin \sqrt{c} t, \frac{2\sqrt{c}}{\sqrt{b^{2}-4}}\sinh s,\sqrt{b^2-4}\cosh s+\! \frac{ab}{\sqrt{{b^{2}\! -\!4}}}\sinh  s\)$, where $c=b^{2}-a^{2}-4 >0$, $f=a\cosh s +b \sinh s$, and $a,b$ are constants;}

\item[{\rm (2)}]  {\it $\phi=\! \frac{1}{\sqrt{c}}\! \(\! \frac{ab}{\sqrt{a^{2}\! + \! 4}}\cosh s\!+\!\sqrt{a^{2}+4}\sinh s,  f \sinh \! \sqrt{c}t,  f\cosh \!\sqrt{c}t, \frac{2\sqrt{c}}{\sqrt{b^2+4}}\cosh s\! \)$, where $c=4+a^{2}-b^{2}>0$, $f=a\cos s +b \sin s$, and $a,b$ are constants;}

\item[{\rm (3)}] {\it $\phi=\(f t, \sinh s-\dfrac{b}{2}f t^{2}, ft^{2},  \cosh s +\dfrac{a}{2}f t\)$,\ where $f=a\cosh s +b \sinh s$, and $a,b$ are constants satisfying $b^{2}-a^{2}=4$;}

\item[{\rm (4)}] {\it  $\phi=\(s\cos t, s\sin t, \sqrt{1+s^2}  \cos f, \sqrt{1+s^2}  \sin f\)$, where 
$$f=\pm \int^{s}\frac{(2+k)ds} {(1+s^2)s^{k/(2+k)}\sqrt{(2+k)^2s^{-2k/(2+k)}+s^2+1}},$$ with
$s\in (0,1)$ and $k\ne 0,-1,-2$;}

\item[{\rm (5)}] {\it  $\phi=\(s\sinh t,  \sqrt{s^2-1}  \cosh f,s\cosh t, \sqrt{s^2-1}  \sinh f\)$, where 
$$f=\pm \int^{s}\frac{(2+k)ds} {(s^2-1)s^{k/(2+k)}\sqrt{(2+k)^{2}s^{-2k/(2+k)}+s^2+1}}$$ with
$s\in (1,\infty)$ and $k\ne 0,-1,-2$;}

\item[{\rm (6)}]  {\it  $\phi=\(s\sinh t, \sqrt{1-s^2}  \sinh f, \sqrt{1-s^2}  \cosh f, s\cosh t\)$, where 
$$f=\pm \int^{s}\frac{(2+k)ds} {(s^2-1)s^{k/(2+k)}\sqrt{(2+k)^2s^{-2k/(2+k)}+s^2+1}},$$ with
$s\in (-1,1)$ and $k\ne 0,-1,-2$;}

\item[{\rm (7)}]  {\it  $\phi=\(\frac{1}{2}\(st^{2}-sf^{2}+\frac{1}{s}-s\), st, sf,\frac{1}{2}\(st^{2}-sf^{2}+\frac{1}{s}+s\)\)$, where 
$$f= \int^{s}\frac{(2+k)ds} {s^2\sqrt{(2+k)^2+s^{(4+4k)/(2+k)}}},$$ with
$s\in (0,\infty)$ and $k\ne 0,-1,-2$.}
\end{itemize}\end{theorem}

\end{document}